\documentclass[letterpaper, preprint, paper,11pt]{AAS}	% for preprint proceedings

\usepackage{bm}
\usepackage{amsmath}
\usepackage{amsfonts}
\usepackage{subfigure}
\usepackage[colorlinks=true, pdfstartview=FitV, linkcolor=black, citecolor= black, urlcolor= black]{hyperref}
\usepackage{overcite}
\usepackage{footnpag}			      	% make footnote symbols restart on each page
\usepackage{algorithm,algcompatible}
\algnewcommand\INPUT{\item[\textbf{Input:}]}%
\algnewcommand\OUTPUT{\item[\textbf{Output:}]}%

\PaperNumber{24-167}

\begin{document}

\title{Time Efficient Rate Feedback Tracking Controller with Slew Rate and Control Constraint}

\author{Seungyeop Han\thanks{PhD Student, Daniel Guggenheim School of Aerospace Engineering, Georgia Institute of Technology, Atlanta, GA,30332},\
Byeong-Un Jo\thanks{Assistant Professor, Department of Aerospace Engineering, Sejong University, Seoul 05006, Republic of Korea},\
 and Koki Ho\thanks{Dutton-Ducoffe Professor, Daniel Guggenheim School of Aerospace Engineering, Georgia Institute of Technology, Atlanta, GA, 30332}
}

\maketitle{}

\begin{abstract}
This paper proposes a time-efficient attitude-tracking controller considering the slew rate constraint and control constraint. The algorithm defines the sliding surface, which is the linear combination of command, body, and regulating angular velocity, and utilizes the sliding surface to derive the control command that guarantees finite time stability. The regulating rate, which is an angular velocity regulating the attitude error between the command and body frame, is defined along the instantaneous eigen-axis between the two frames to minimize the rotation angle. In addition, the regulating rate is shaped such that the slew rate constraint is satisfied while the time to regulation is minimized with consideration of the control constraint. Practical scenarios involving Earth observation satellites are used to validate the algorithm's performance.
\end{abstract}

\section*{Notation}
\begin{tabular}{r l}
	$\textbf{\textit{a}}$  & geometric vector \\
    $\hat{\textbf{\textit{a}}}$ & unit vector of $\textbf{\textit{a}}$ \\
    $\textit{\textit{a}}$ & magnitude of $\textbf{\textit{a}}$ or scalar value \\
    $\textbf{\textit{a}}^A$ & $\textbf{\textit{a}}$ expressed in frame $A$ \\
    $\dot{a}$ & time derivative of scalar $a$ \\
    $^{A} \dot{\textbf{\textit{a}}}$ & time derivative of $\textbf{\textit{a}}$ with respect to frame $A$ \\
    $^{A} \dot{\textbf{\textit{a}}}^A$ & $^{A} \dot{\textbf{\textit{a}}}$ expressed in frame $A$ \\
\end{tabular} \\

\section{Introduction}
The faster maneuvering capability of Earth observation satellites (EOS) has always been an important research area. The satellite's agility directly indicates its operational capacity per orbit pass, and the capability of agile tracking maneuver allows EOS to have various imaging operations, such as fixed point staring imaging and non-parallel ground scan imaging. Moreover, modern satellite systems often necessitate autonomous and agile responsive operational capabilities. Therefore, the rotational agility facilitated by onboard attitude control algorithms serves as a crucial performance factor in the development of modern satellite systems.

Numerous studies have been conducted on this topic, and it remains an active area of research today. In the early days, well-known quaternion-based attitude control laws were proposed as non-singular solutions, either in the form of sliding mode control (SMC) or PD control \cite{vadali1986variable, wie1989quarternion}. Soon after, attitude-tracking control laws were also developed\cite{weiss1993quaternion, crassidis2000optimal}. Despite their simple structure, these methods were very successful, leading to subsequent studies that built upon and slightly modified them to achieve more desirable properties. For example, some papers considered slew rate constraints \cite{wie1995feedback, verbin2011time, cao2016time, schaub2018speed}, some studies achieved finite-time stability \cite{jin2008robust, wu2011quaternion, li2017robust}, and others considered model parameters to have time-efficient maneuver characteristics \cite{wie2002rapid, verbin2011time, cao2016time, li2018time}.

Most feedback-type control laws adopt an eigen-axis rotation structure for two main reasons. First, it is known that attitude rotation along the eigen-axis is time-efficient and close to time-optimal for rest-to-rest maneuvers if the rotation angle is not significant \cite{bilimoria1993time}. Second, the rotational motion along the eigen-axis can be expressed analytically. The optimal control profile for time-optimal problems has a Bang-Bang structure \cite{kirk2004optimal}, and using the eigen-axis enables obtaining an explicit solution under the Bang-Bang structure. Some works utilize these properties to generate attitude trajectories in a closed-form manner \cite{schaub1996globally, you2019near}.

This paper aims to develop a time-efficient attitude-tracking control law that considers both slew rate and acceleration constraints. Like many previous studies, it adopts a finite-time convergence SMC structure. However, the structure is rearranged to represent the desired body frame rate explicitly. This modification allows for shaping a time-efficient rate profile while easily handling the constraints. The paper adopts a trapezoidal acceleration profile to make the maneuver time-efficient, which alleviates the infinite jerk issue and guarantees finite-time convergence. However, this profile does not work well for low-frequency control, so this paper additionally introduces a modified trapezoidal profile that performs well in practical systems.

The remainder of the paper is structured as follows. First, we review the problem statement and the kinematics of the error quaternion. Next, we propose the structure of the sliding surface and adopt a well-known SMC method to achieve finite-time convergence toward the sliding manifold. Then, we elucidate the design of the regulating rate profile, which can handle both slew rate and acceleration constraints. Based on the previously explained concepts, we then outline the overall sequence for computing the time-efficient torque command. The effectiveness of the proposed method is demonstrated through numerical simulations of various practical scenarios. Finally, we conclude the paper with a discussion of future work.

%----------------------------------------------------------------------------------------
%	Rate Feedback Controller
%----------------------------------------------------------------------------------------
\section{Rate Feedback Controller}
%----------------------------------------------------------------------------------------
%	Problem Statement
%----------------------------------------------------------------------------------------
\subsection{Problem Statement}
This paper uses a unit quaternion to represent an attitude. The attitude of frame $B$ with respect to frame $A$ can be expressed as:
\begin{equation}
    \bar{\textbf{\textit{q}}}_{B/A} = \begin{bmatrix}
        \textbf{\textit{q}}_{B/A} \\ q_{B/A}
    \end{bmatrix},\quad \lVert \bar{\textbf{\textit{q}}}_{B/A} \rVert = 1
\end{equation}
Based on the quaternion definition we used, the attitude kinematics is given as:
\begin{equation}
    \dot{\bar{\textbf{\textit{q}}}}_{B/A}= \frac{1}{2} \bar{\boldsymbol{\omega}}_{B/A} \otimes 
    \bar{\textbf{\textit{q}}}_{B/A} = \frac{1}{2} \begin{bmatrix}
        \textit{q}_{B/A} \boldsymbol{\omega}_{B/A}  - \boldsymbol{\omega}_{B/A} \times \textbf{\textit{q}}_{B/A}  \\
        - \boldsymbol{\omega}_{B/A}  \cdot \textbf{\textit{q}}_{B/A} 
    \end{bmatrix}
\end{equation}
where $\bar{\boldsymbol{\omega}}_{B/A} = [{\boldsymbol{\omega}}_{B/A}^\top\ 0]^\top$ and full expression of ${\boldsymbol{\omega}}_{B/A}$ of the equation is ${\boldsymbol{\omega}}_{B/A}^B$. In addition, $\otimes$ is quaternion operator defined for two quaternions $\bar{p}$ and $\bar{q}$ as:
\begin{equation}
    \bar{\textbf{\textit{p}}} \otimes \bar{\textbf{\textit{q}}} = \begin{bmatrix}
        p \textbf{\textit{q}} + q \textbf{\textit{p}} - \textbf{\textit{p}} \times \textbf{\textit{q}}\\
        pq - \textbf{\textit{p}}\cdot \textbf{\textit{q}}
    \end{bmatrix}
\end{equation}
Lastly, if $\mathcal{I}$ is used as the base frame then we will omit the expression for simplicity, i.e. $\bar{\textbf{\textit{q}}}_{\mathcal{B}/\mathcal{I}} = \bar{\textbf{\textit{q}}}_{\mathcal{B}}$ and $
\boldsymbol{\omega}_{\mathcal{B}/\mathcal{I}} = \boldsymbol{\omega}_{\mathcal{B}}$.

The spacecraft attitude dynamics is given as
\begin{equation}
    J \dot{\boldsymbol{\omega}}_\mathcal{B} = \textbf{\textit{u}} + \textbf{\textit{d}} - \boldsymbol{\omega}_\mathcal{B} \times 
    J \boldsymbol{\omega}_\mathcal{B}
\end{equation}
where $J$ is the inertia matrix of spacecraft in $\mathcal{B}$ frame, $\textbf{\textit{u}}$ and $\textbf{\textit{d}}$ represent the control command and external disturbance torque, respectively. We assume that the norm of $\textbf{\textit{d}}$ is bounded by $d_\text{max}$ and known. Note that full expression of $\dot{\boldsymbol{\omega}}_\mathcal{B}$ is ${}^\mathcal{B}\dot{\boldsymbol{\omega}}_\mathcal{B}^\mathcal{B}$, i.e., the time derivative of $\boldsymbol{\omega}$ in $\mathcal{B}$ and expressed in $\mathcal{B}$.

Let $\mathcal{D}$ be the command/desired frame to be tracked by $\mathcal{B}$, and the kinematics of $\mathcal{D}$ is given as:
\begin{equation}
    \dot{\bar{\textbf{\textit{q}}}}_\mathcal{D}= \frac{1}{2} \bar{\boldsymbol{\omega}}_\mathcal{D} \otimes 
    \bar{\textbf{\textit{q}}}_\mathcal{D}
\end{equation}
The goal of the paper is to derive $\textbf{\textit{u}}$ that make $
\bar{\textbf{\textit{q}}}_\mathcal{B} \to \bar{\textbf{\textit{q}}}_\mathcal{D}
$ and $\boldsymbol{\omega}_\mathcal{B} \to \boldsymbol{\omega}_\mathcal{D}$ as $t\to T < \infty$ for some $T>0$ while ensuring $\lVert \boldsymbol{\omega}_\mathcal{B}(t) \rVert_2 < \omega_\text{max}$ and $\lVert \textbf{\textit{u}}\rVert_2 < u_\text{max}$ for all $t$.

%----------------------------------------------------------------------------------------
%	Error Quaternion and Kinematics
%----------------------------------------------------------------------------------------
\subsection{Error Quaternion and Kinematics}
We define error quaternion and angular velocity as follows:
\begin{equation} \label{eq:error_state}
    \bar{\textbf{\textit{q}}}_e \equiv \bar{\textbf{\textit{q}}}_\mathcal{D/B} = 
    \bar{\textbf{\textit{q}}}_\mathcal{D} \otimes
    \bar{\textbf{\textit{q}}}_\mathcal{B}^{-1}, \quad
    \boldsymbol{\omega}_e \equiv
    \boldsymbol{\omega}_\mathcal{D} - \boldsymbol{\omega}_\mathcal{B}
\end{equation}
where $\bar{\textbf{\textit{q}}}^{-1} = [-\textbf{\textit{q}}^\top\ q]^\top$. Note that it is better to express a vector in $\mathcal{B}$ for later usage. For example, if $\omega_\mathcal{D}^\mathcal{D}$ is given then $\omega_e^\mathcal{B} = T_{\mathcal{B}/\mathcal{D}}\omega_\mathcal{D}^\mathcal{D}- \omega_\mathcal{B}^\mathcal{B}$ where $T_{A/B}$ is the frame rotational matrix from $B$ to $A$. Then the time derivative of error quaternion can be computed as:
\begin{equation} \label{eq:error_kinematics}
\begin{aligned}\dot{\bar{\textbf{\textit{q}}}}_e &= 
\dot{\bar{\textbf{\textit{q}}}}_\mathcal{D} \otimes \bar{\textbf{\textit{q}}}_\mathcal{B}^{-1} + \bar{\textbf{\textit{q}}}_\mathcal{D} \otimes \dot{\bar{\textbf{\textit{q}}}}_\mathcal{B}^{-1}\\
&=\frac{1}{2} \bar{\boldsymbol{\omega}}_\mathcal{D}^\mathcal{D} \otimes 
    \bar{\textbf{\textit{q}}}_\mathcal{D} \otimes \bar{\textbf{\textit{q}}}_\mathcal{B}^{-1} 
    - \frac{1}{2} \bar{\textbf{\textit{q}}}_\mathcal{D} 
    \otimes \bar{\textbf{\textit{q}}}_\mathcal{B}^{-1} \otimes \bar{\boldsymbol{\omega}}_\mathcal{B}^\mathcal{B} \\
    &= \frac{1}{2} \bar{\boldsymbol{\omega}}_\mathcal{D}^\mathcal{D} \otimes 
    \bar{\textbf{\textit{q}}}_e
    - \frac{1}{2} \bar{\textbf{\textit{q}}}_\mathcal{D/B} \otimes \bar{\boldsymbol{\omega}}_\mathcal{B}^\mathcal{B} \otimes \bar{\textbf{\textit{q}}}_\mathcal{D/B}^{-1} \otimes \bar{\textbf{\textit{q}}}_\mathcal{D/B}\\
    &= \frac{1}{2} \left( \bar{\boldsymbol{\omega}}_\mathcal{D}^\mathcal{D} - \bar{\boldsymbol{\omega}}_\mathcal{B}^\mathcal{D} \right) \otimes \bar{\textbf{\textit{q}}}_e \\
    &= \frac{1}{2} \bar{\boldsymbol{\omega}}_e^\mathcal{D} \otimes \bar{\textbf{\textit{q}}}_e
\end{aligned}
\end{equation}
Note that the instantaneous eigen-axis and the attitude angle error can be computed as follows:
\begin{equation} \label{eq:eigenaxis}
    \hat{\textbf{\textit{e}}}^\mathcal{D} = \hat{\textbf{\textit{e}}}^\mathcal{B} = \frac{\textbf{\textit{q}}_e}{\lVert \textbf{\textit{q}}_e \rVert},\quad \theta_e = 2\cos^{-1} (q_e)
\end{equation}
and the error quaternion can also be expressed as follows:
\begin{equation} 
    \bar{\textbf{\textit{q}}}_e = \begin{bmatrix}
        {\textbf{\textit{q}}}_e \\ q_e
    \end{bmatrix} =
    \begin{bmatrix}
        \sin\left(\frac{\theta_e}{2}\right) \hat{\textbf{\textit{e}}}^\mathcal{B}  \\
        \cos\left(\frac{\theta_e}{2}\right) 
    \end{bmatrix}
\end{equation}
Lastly, the time derivative of the error angle is:
\begin{equation} \label{eq:errorangle_dot}
    \dot \theta_e = -\frac{2 }{\sqrt{1-q_e^2}} \dot q_e =  \boldsymbol{\omega}_e \cdot \frac{\textbf{\textit{q}}_e}{\lVert \textbf{\textit{q}}_e \rVert} =
    \boldsymbol{\omega}_e  \cdot \hat{\textbf{\textit{e}}}
\end{equation}
and the time derivative of the eigen-axis in frame $\mathcal{B}$ is:
\begin{equation} \label{eq:eigenaxis_dot}
\begin{aligned}
    \dot {\hat{\textbf{\textit{e}}}} &= 
    \left(
    \textbf{\textit{q}}_e \cdot {\textbf{\textit{q}}}_e \right)
    \frac{\dot{\textbf{\textit{q}}}_e}{\lVert \textbf{\textit{q}}_e \rVert^3}
    - \left(
    \textbf{\textit{q}}_e \cdot \dot{\textbf{\textit{q}}}_e \right) \frac{\textbf{\textit{q}}_e}{\lVert \textbf{\textit{q}}_e \rVert^3}
    \\
    &= \frac{1}{\lVert \textbf{\textit{q}}_e \rVert^3} \textbf{\textit{q}}_e \times \left( \dot{\textbf{\textit{q}}}_e \times \textbf{\textit{q}}_e \right) \\
    &= \frac{1}{2\lVert \textbf{\textit{q}}_e \rVert^3} \textbf{\textit{q}}_e \times \left(  q_e \boldsymbol{\omega}_e^\mathcal{B} \times \textbf{\textit{q}}_e - (\textbf{\textit{q}}_e \cdot \textbf{\textit{q}}_e) \boldsymbol{\omega}_e^\mathcal{B} + (\boldsymbol{\omega}_e \cdot \textbf{\textit{q}}_e) \textbf{\textit{q}}_e  \right) \\
    &= \frac{1}{2} \left( \frac{q_e}{\lVert \textbf{\textit{q}}_e \rVert} \boldsymbol{\omega}_{e_\perp} + \boldsymbol{\omega}_{e_\perp} \times \hat{\textbf{\textit{e}}}\right)
\end{aligned}
\end{equation}
where $\boldsymbol{\omega}_{e_\perp} \equiv \boldsymbol{\omega}_e - 
(\hat{\textbf{\textit{e}}} \cdot \boldsymbol{\omega}_e )\boldsymbol{\omega}_e $ is perpendicular decomposition of $\boldsymbol{\omega}_{e}$ with respect to $\hat{\textbf{\textit{e}}}$. Note that the vector triple product identity $\textbf{\textit{a}}\times(\textbf{\textit{b}} \times \textbf{\textit{c}})= (\textbf{\textit{a}}\cdot \textbf{\textit{c}}) \textbf{\textit{b}} - (\textbf{\textit{a}}\cdot \textbf{\textit{b}})\textbf{\textit{c}}$ is used in the first line of the equation.

\textbf{Remark 1.} If $\theta_e = 0$ then $\hat{\textbf{\textit{e}}}$ becomes singular. Additionally, as $\theta_e \to 0$, $\dot {\hat{\textbf{\textit{e}}}}$ becomes numerically unstable and can blow up if $\boldsymbol{\omega}_e \to 0$ not fast enough. If a control law uses either $\hat{\textbf{\textit{e}}}$ or $\dot {\hat{\textbf{\textit{e}}}}$ terms, proper analysis is required, and these issues will be addressed in a later section.

\textbf{Remark 2.} If $\boldsymbol{\omega}_{e_\perp} = 0$, then $\dot {\hat{\textbf{\textit{e}}}} = 0 $ even if $\boldsymbol{\omega}_{e} \neq 0$. That is, if the off-axis error rate $\boldsymbol{\omega}_{e_\perp}$ is 0, then the eigen-axis becomes constant, and the relative motion between $\mathcal{D}$ and $\mathcal{B}$ becomes pure single-axis rotational motion.

%----------------------------------------------------------------------------------------
%	Sliding Surface and Its Stability
%----------------------------------------------------------------------------------------
\subsection{Sliding Surface and Its Stability}
Let the sliding surface $\textbf{\textit{s}}$ be:
\begin{equation}
    \textbf{\textit{s}} = \boldsymbol{\omega}_\mathcal{D} + \boldsymbol{\omega}_R - \boldsymbol{\omega}_\mathcal{B}
\end{equation}
where $\boldsymbol{\omega}_R$ is the regulating rate which will be defined later. Then, the time derivative of the sliding vector in frame $\mathcal{B}$ is:
\begin{equation}
\begin{aligned}
    {}^\mathcal{B} \dot{\textbf{\textit{s}}} &= {}^\mathcal{B} \dot{\boldsymbol{\omega}}_\mathcal{D} + {}^\mathcal{B}\dot{\boldsymbol{\omega}}_R - {}^\mathcal{B}\dot{\boldsymbol{\omega}}_\mathcal{B} \\
    &= {}^\mathcal{B} \dot{\boldsymbol{\omega}}_\mathcal{D} + {}^\mathcal{B}\dot{\boldsymbol{\omega}}_R - J^{-1}\left( \textbf{\textit{u}} + \textbf{\textit{d}} - \boldsymbol{\omega}_\mathcal{B} \times 
    J \boldsymbol{\omega}_\mathcal{B} \right)
\end{aligned}
\end{equation}
Define the control command $\textbf{\textit{u}}$ as
\begin{equation}\label{eq:control_law} 
    \textbf{\textit{u}} = J\left( {}^\mathcal{B} \dot{\boldsymbol{\omega}}_\mathcal{D} + {}^\mathcal{B}\dot{\boldsymbol{\omega}}_R +  \beta_1 s^{\beta_2} \hat{\textbf{\textit{s}}}  \right) + d_\text{max} \hat{\textbf{\textit{s}}} + \boldsymbol{\omega}_\mathcal{B} \times 
    J \boldsymbol{\omega}_\mathcal{B}
\end{equation}
where $\beta_1 > 0$ and $0 < \beta_2 < 1$ are constant control parameters to be tuned by a user. The unit vector $\hat{\textbf{\textit{s}}}$ should be computed as follows:
\begin{equation}
    \hat{\textbf{\textit{s}}} = \begin{cases}
        \frac{\textbf{\textit{s}}}{s} & \mbox{if } s > 0 \\
        0 & \mbox{if } s = 0 \\
    \end{cases}
\end{equation}
to avoid the division by zero.

\textbf{Proposition 1.} The control law of Eq.~\eqref{eq:control_law} dervies the states toward sliding manifold $\textbf{\textit{s}}= 0$ within finite time.

\noindent \textbf{Proof.} Define the valid Lyapunov function as
\begin{equation}
    V = \frac{1}{2}\textbf{\textit{s}}^\top J\textbf{\textit{s}}
\end{equation}
then the time derivative of $V$ is
\begin{equation}
\begin{aligned}
    \dot V &= \textbf{\textit{s}}^\top J\dot{\textbf{\textit{s}}} \\
    &= \textbf{\textit{s}}^\top J\left({}^\mathcal{B} \dot{\boldsymbol{\omega}}_\mathcal{D} + {}^\mathcal{B}\dot{\boldsymbol{\omega}}_R - J^{-1}\left( \textbf{\textit{u}} + \textbf{\textit{d}} - \boldsymbol{\omega}_\mathcal{B} \times 
    J \boldsymbol{\omega}_\mathcal{B} \right)  \right) \\
    &= -\textbf{\textit{s}}^\top J\left( J^{-1}(  d_\text{max} \hat{\textbf{\textit{s}}} + \textbf{\textit{d}} ) + \beta_1 s^{\beta_2} \hat{\textbf{\textit{s}}} 
     \right) \\
    &< -\beta_1 \lambda_\text{min}(J) s^{\beta_2 + 1} = - \beta_1 \lambda_\text{min}(J) V^{\frac{\beta_2+1}{2}}
\end{aligned}
\end{equation}
and this guarantees finite time stability toward the sliding surface since $\beta_2 < 1$. 

\textbf{Remark 3.} If states are on the manifold, the error quaternion is directly controlled by the $\boldsymbol{\omega}_R$ since 
\begin{equation}
    \textbf{\textit{s}} = \boldsymbol{\omega}_\mathcal{D} + \boldsymbol{\omega}_R - \boldsymbol{\omega}_\mathcal{B} = 0 \implies 
    \boldsymbol{\omega}_e = 
    \boldsymbol{\omega}_\mathcal{D} - \boldsymbol{\omega}_\mathcal{B} = -\boldsymbol{\omega}_R
\end{equation}
and by the Eq.~\eqref{eq:error_kinematics}.
By properly designing the $\boldsymbol{\omega}_R$, the error states can be regulated time efficiently while considering rate and control limits.

\textbf{Remark 4.} Although the control law is robust against disturbances, it inherently has a chattering problem. This issue can be mitigated by adopting an improved sliding mode structure or having an independent disturbance observer, which will be discussed in future work.

% Reconsider this for journal paper
% \textbf{Remark 5.} $\hat{\textbf{\textit{s}}}$ becomes singular when $\textbf{\textit{s}} = 0$. In practice, $\textbf{\textit{s}}$ will oscillate around 0 due to sensor noise from $\boldsymbol{\omega}_\mathcal{B}$. Therefore, we define the normal vector as 
% \begin{equation}
%     \hat{\textbf{\textit{s}}} = \frac{\textbf{\textit{s}}}{\max(s,\varepsilon)}
% \end{equation}
% where $\varepsilon > 0$ is a small value, approximately equal to the standard deviation of the rate sensor noise.

%----------------------------------------------------------------------------------------
%	Design of Regulating Rate and Its Properties
%----------------------------------------------------------------------------------------
\subsection{Design of Regulating Rate and Its Properties}
This work also adopts rotation along the eigen-axis to achieve time-efficient maneuvers, as seen in many previous studies. Additionally, if the regulating angular velocity is designed to be aligned with the eigen-axis $\hat{\textbf{\textit{e}}}$, then the relative motion becomes a pure single-axis rotational problem, allowing for a closed-form expression. For these reasons, we define
\begin{equation}
    \boldsymbol{\omega}_R = \omega_R \hat {\textbf{\textit{e}}}
\end{equation}
where $\omega_R$ is the regulating rate to be defined.

\subsubsection{Maximum Acceleration Level} 
The $\omega_R$ profile must satisfy the control constraint, i.e., the acceleration required to track the profile must be admissible with some margin. Based on the control law, the necessary feed-forward accelerations along the curve are the angular acceleration of $\mathcal{D}$ and the gyroscopic acceleration. Therefore, the maximum acceleration along the eigen-axis is defined as follows:
\begin{equation} \label{eq:alpha_max}
    \alpha_\text{max}(\hat{\textbf{\textit{e}}}) = \frac{\gamma}{\lVert J \hat{\textbf{\textit{e}}} \rVert}\left( u_\text{max} - \lVert J {}^\mathcal{B} \dot{\boldsymbol{\omega}}_\mathcal{D} \rVert - \lVert \boldsymbol{\omega}_\mathcal{B} \times J \boldsymbol{\omega}_\mathcal{B} \rVert \right)
\end{equation}
where $\gamma \approx 1$ is the user-defined deceleration ratio that considers the uncertainty of system parameters.

Note that $\alpha_\text{max}$ is a function of $\hat{\textbf{\textit{e}}}$, which may cause problems when $\theta_e$ is very small because small perturbation on $\textbf{\textit{q}}_e$ can greatly change the eigen-axis and, consequently, $\alpha_\text{max}$ as well. For this reason, we define the regulating acceleration as follows:
\begin{equation} \label{eq:alpha_R}
    \alpha_R = \left( 1 - \sigma(\theta_e; \eta) \right) \alpha_\text{min} + \sigma(\theta_e; \eta) \alpha_\text{max}
\end{equation}
where $\alpha_\text{min}$ is the minimum value of $\alpha_\text{max}$, which is the maximum acceleration along the major axis, defined as:
\begin{equation}
    \alpha_\text{min} = \min_{\hat{\textbf{\textit{x}}}} \alpha_\text{max} (\hat{\textbf{\textit{x}}})
\end{equation}
where $\hat{\textbf{\textit{x}}}$ is a general unit vertor. In addition, $\sigma(\cdot,\eta)$ is a linear saturation function with a parameter $\eta > 0$ defined as:
\begin{equation}
    \sigma(x; \eta) = \begin{cases}
        \frac{x}{\eta} & \mbox{if } x < \eta \\
        1 & \mbox{if } x \geq \eta \\
    \end{cases}
\end{equation}
for $x>0$. This modification prevents fluctuations in the $\alpha_R$ value when $\theta_e$ is small.

\subsubsection{Maximum Regulating Rate} In order to ensure $\omega_\mathcal{B} \leq \omega_\text{max}$, the maximum of $\omega_R$ must satisfy $\lVert \boldsymbol{\omega}_\mathcal{D} + \boldsymbol{\omega}_R \rVert \leq \omega_\text{max}$. Let $\omega_{R_\text{max}}$ be the maximum value satisfying the equality, then it can be solved as:
\begin{equation} \label{eq:omega_max}
    \omega_{R_\text{max}} = -\left( \boldsymbol{\omega}_\mathcal{D} \cdot \hat{\textbf{\textit{e}}} \right) + \sqrt{\left( \boldsymbol{\omega}_\mathcal{D} \cdot \hat{\textbf{\textit{e}}} \right)^2 + \omega_\text{max}^2 - \left( \boldsymbol{\omega}_\mathcal{D} \cdot \boldsymbol{\omega}_\mathcal{D} \right)}
\end{equation}
with assumption of $\lVert \boldsymbol{\omega}_\mathcal{D} \rVert < \omega_\text{max}$. Note that the time derivative of $\omega_{R_\text{max}}$ is
\begin{equation} \label{eq:omega_max_dot}
    \dot \omega_{R_\text{max}} =  -\dot{\omega}_{\mathcal{D}_e} + \frac{
    \left( \boldsymbol{\omega}_\mathcal{D} \cdot \hat{\textbf{\textit{e}}} \right)\dot{\omega}_{\mathcal{D}_e} - \dot{\boldsymbol{\omega}}_\mathcal{D} \cdot \boldsymbol{\omega}_\mathcal{D}}{\sqrt{\left( \boldsymbol{\omega}_\mathcal{D} \cdot \hat{\textbf{\textit{e}}} \right)^2 + \omega_\text{max}^2 - \left( \boldsymbol{\omega}_\mathcal{D} \cdot \boldsymbol{\omega}_\mathcal{D} \right)}}
\end{equation}
where $\dot{\omega}_{\mathcal{D}_e} =  \dot{\boldsymbol{\omega}}_\mathcal{D} \cdot \hat{\textbf{\textit{e}}} + {\boldsymbol{\omega}}_\mathcal{D} \cdot \dot{\hat{\textbf{\textit{e}}}}$. It is clear that $ \omega_{R_\text{max}} = \omega_\text{max}$ and $\dot \omega_{R_\text{max}} = 0$ when the desired frame is inertially fixed.

% 저널버젼에서는 Wd_dot이 무시가능하고, Wb(t_0) < wmax 였으면, Wb(t) < wmax임을 보이면 좋을 듯. Wb(t) = Wmax인 경우 Wb_dot * Wb <= 0 임을 보여서 적분 관계 보이면 될듯, 추가로 더 조건이 필요할듯
% \textbf{Proposition 2.} The control law of Eq.~\eqref{eq:control_law} dervies the states toward sliding manifold $\textbf{\textit{s}}= 0$ within finite time.

\subsubsection{Time-Efficient Rate Profile}
When the acceleration is bounded, the time-optimal control profile has a Bang-Bang structure[]. However, such a profile is sensitive to any uncertainty and requires an infinite level of jerk. Therefore, a trapezoidal profile is preferred, with little sacrifice in time-optimality. However, directly using a trapezoidal profile for $\omega_R$ will also cause practical problems. In this paper, we introduce a trapezoidal profile and a modified trapezoidal profile.

\begin{figure}[h]
    \centering
    \includegraphics[width=0.4\textwidth]{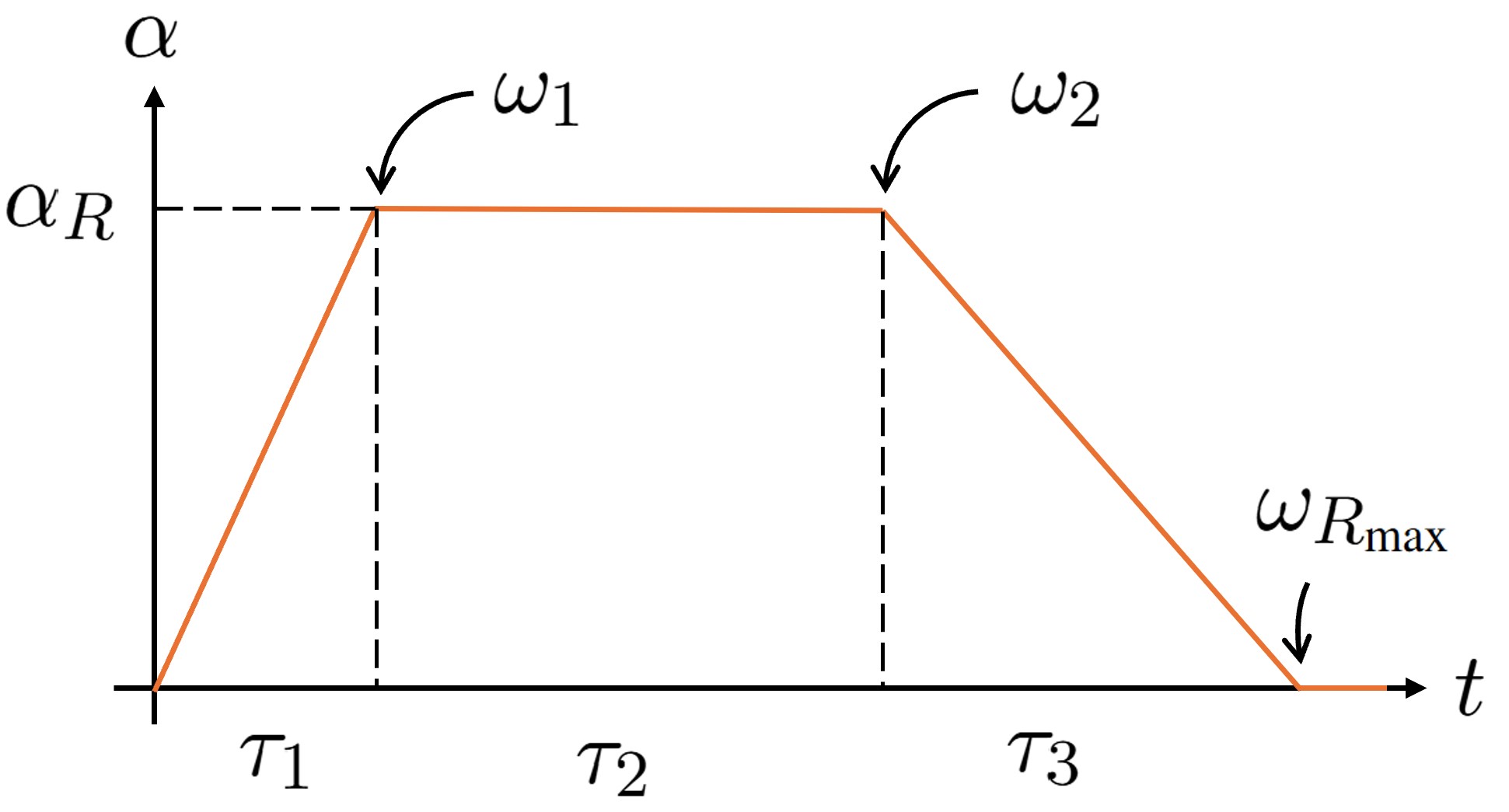}
    \includegraphics[width=0.4\textwidth]{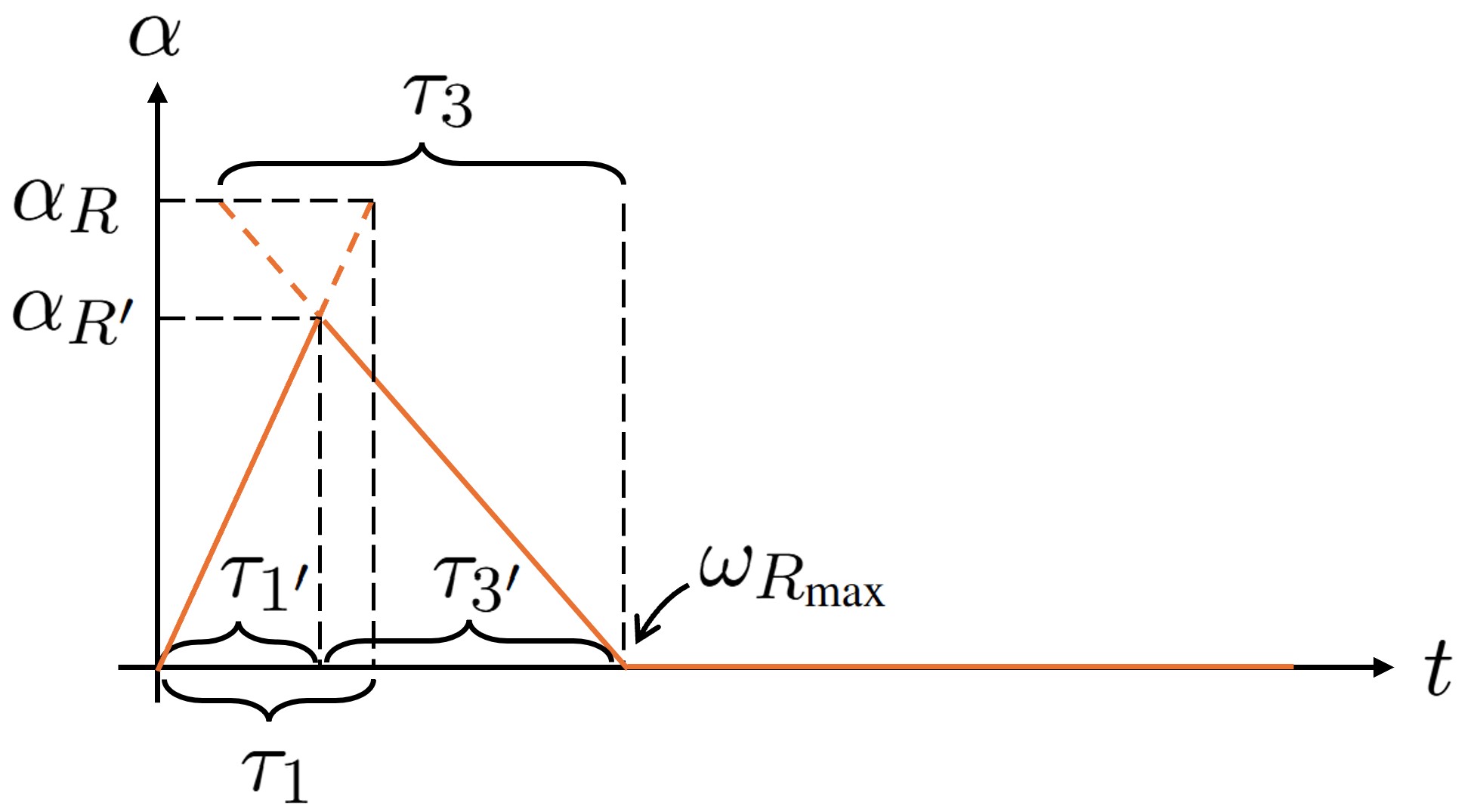}
    \caption{Trapezoidal acceleration profile for (a) $\tau_2>0$ (b) $\tau_2 < 0$}
    \label{fig:trape_profile}
\end{figure}

\textbf{Trapezoidal Profile: }
Let $\tau_1$ and $\tau_3$ be the durations for the linearly varying acceleration segments, which are user-defined parameters. We know that the rate is 0 at the beginning of the first segment and becomes $\omega_{R_\text{max}}$ at the end of the third segment, as shown in Figure~\ref{fig:trape_profile}(a). We first compute the rate at the end of the first and second segments as:
\begin{equation}
    \omega_{1} = \frac{1}{2} \alpha_R \tau_1, \quad \omega_{2} = \omega_{R_\text{max}} - \frac{1}{2}\alpha_R\tau_3
\end{equation}
Then the duration of the second segment, which has constant plateau acceleration, is 
\begin{equation}
    \tau_2 = \frac{\omega_{2} - \omega_{1}}{\alpha_R}
\end{equation}
Note that depending on the value of the parameters, $\tau_2$ can be negative, and that is the case where there is no plateau segment as shown in Figure.~\ref{fig:trape_profile}(b). 

If $\tau_2 > 0$, compute the $\omega_R$ as follows:
\begin{equation} \label{eq:wr_profile}
    \omega_R = \begin{cases}
        \frac{1}{2}\frac{\alpha_R}{\tau_1}\tau^2, \quad \tau = \left(\frac{6\theta_e \tau_1}{a_R}\right)^{\frac{1}{3}} & \mbox{if } 0 \leq \theta_e < \theta_1 \\
        \omega_1 + \alpha_R \tau, \quad \tau = \frac{-\omega_1 + \sqrt{\omega_1^2 + 2\alpha_R(\theta_e - \theta_1)}}{\alpha_R} & \mbox{if } \theta_1 \leq  \theta_e < \theta_2 \\
        \omega_{R_\text{max}} - \frac{1}{2}\frac{\alpha_R}{\tau_3}\tau^2, \quad \tau = 2\sqrt{-\frac{p}{3}}\cos\left(\frac{1}{3}\cos^{-1}\left(\frac{3q}{2p} \sqrt{\frac{-3}{p}}\right) -\frac{2\pi}{3}\right) & \mbox{if } \theta_2 \leq \theta_e < \theta_3 \\ \omega_{R_\text{max}} & \mbox{else}
    \end{cases}
\end{equation}
where $\theta_1$, $\theta_2$, and $\theta_3$ are computed as 
\begin{equation}
    \theta_{1} = \frac{1}{6}\alpha_R\tau_1^2,\quad \theta_{2} = \theta_{1} + \omega_{1}\tau_{2} + \frac{1}{2}\alpha_R \tau_2^2, \quad \theta_3 = \theta_2 + \omega_2\tau_3 + \frac{1}{3}\alpha_R \tau_3^2
\end{equation}
, and $p$, $q$ for the region $\theta_2 \leq \theta_e < \theta_3$ are
\begin{equation}
    p = -\frac{6\omega_{R_\text{max}}\tau_3}{\alpha_R},\quad q = \frac{6\tau_3(\theta_3 - \theta_e) }{\alpha_R}
\end{equation} 
The detailed derivation is explained in the Appendix.

If $\tau_2 < 0$ then $\alpha_{R'}$, $\tau_{1'}$, and $\tau_{3'}$ are first computed as:
\begin{equation}
    \alpha_{R'} = \sqrt{\frac{2 \alpha_R \omega_{R_\text{max}}}{\tau_1 + \tau_3}},\quad
    \tau_{1'} = \frac{\alpha_{R'}}{\alpha_R}\tau_1, \quad
    \tau_{3'} = \frac{\alpha_{R'}}{\alpha_R}\tau_3
\end{equation}
Integrating the profile gives the $\omega_R$ profile as follows:
\begin{equation}
    \omega_R = \begin{cases}
        \frac{1}{2}\frac{\alpha_{R'}}{\tau_{1'}}\tau^2, \quad \tau = \left(\frac{6\theta_e \tau_{1'}}{a_{R'}}\right)^{\frac{1}{3}} & \mbox{if } 0 \leq \theta_e < \theta_1 \\
        \omega_{R_\text{max}} - \frac{1}{2}\frac{\alpha_{R'}}{\tau_{3'}}\tau^2, \quad \tau = 2\sqrt{-\frac{p}{3}}\cos\left(\frac{1}{3}\cos^{-1}\left(\frac{3q}{2p} \sqrt{\frac{-3}{p}}\right) -\frac{2\pi}{3}\right) & \mbox{if } \theta_1 \leq \theta_e < \theta_2 \\ \omega_{R_\text{max}} & \mbox{else}
    \end{cases}
\end{equation}
where $\theta_1$, and $\theta_2$ are computed as 
\begin{equation}
    \theta_{1} = \frac{1}{6}\alpha_{R'}\tau_{1'}^2,\quad
    \omega_{1} = \frac{1}{2}\alpha_{R'}\tau_{1'},\quad
    \theta_{2} = \theta_{1} + \omega_{1}\tau_{3'} + \frac{1}{3}\alpha_{R'} \tau_{3'}^2
\end{equation}
,and $p$, $q$ for the region $\theta_1 \leq \theta_e < \theta_2$ are
\begin{equation}
    p = -\frac{6\omega_{R_\text{max}}\tau_{3'}}{\alpha_{R'}},\quad q = \frac{6\tau_{3'}(\theta_2 - \theta_e) }{\alpha_{R'}}
\end{equation}
Although the derivation is not presented, a slight modification of the derivation of  Eq.~\eqref{eq:wr_profile} will yield the rate profile for $\tau_2 < 0$.

\begin{figure}[h]
    \centering
    \includegraphics[width=0.45\textwidth]{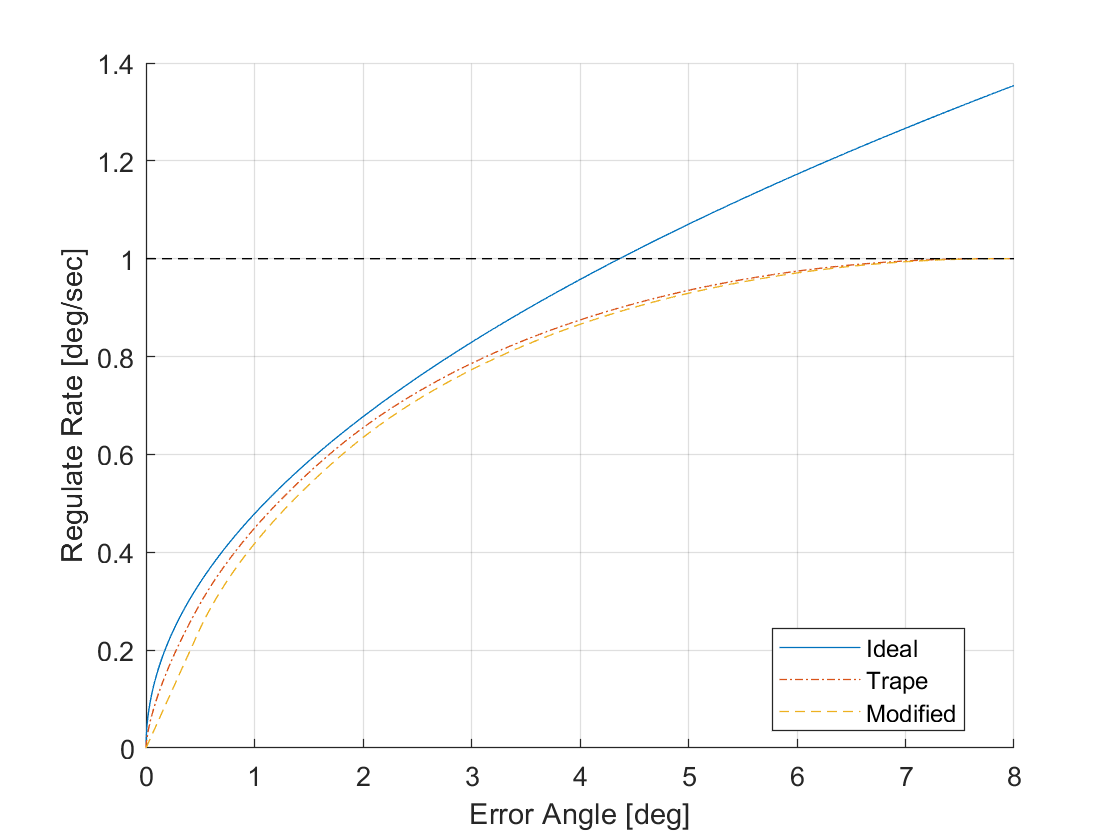}
    \includegraphics[width=0.45\textwidth]{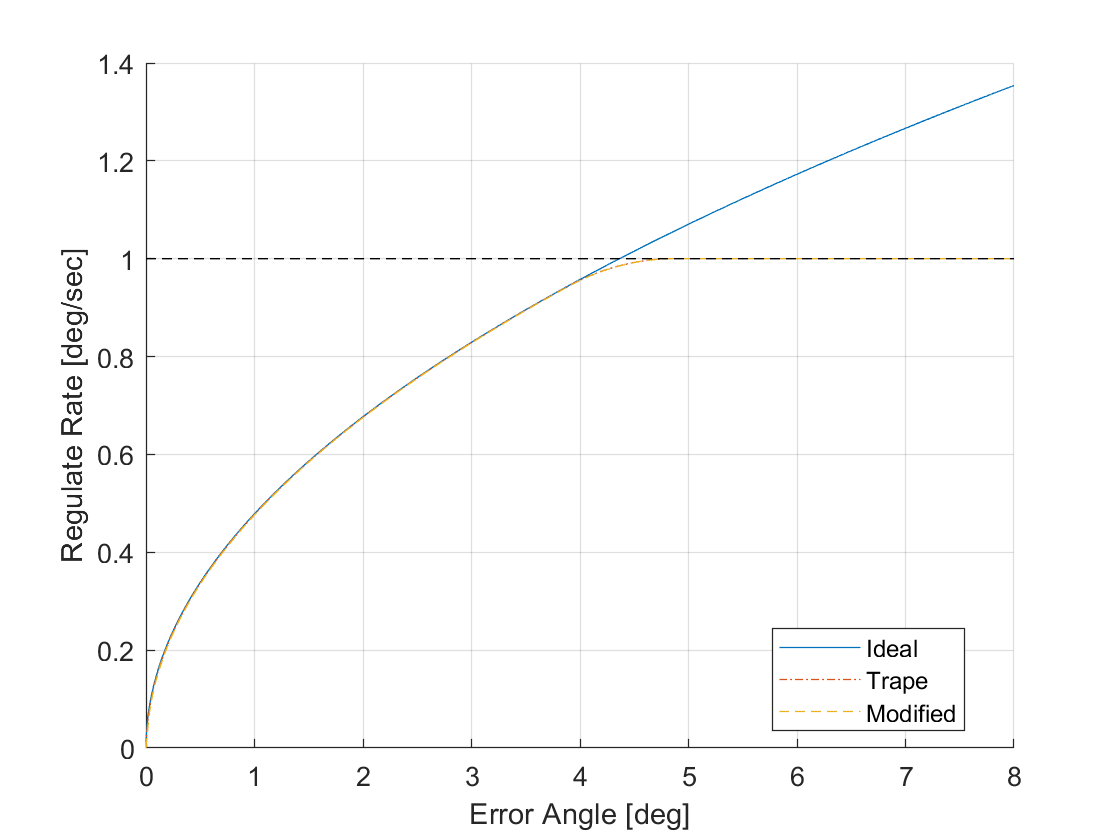}
    \caption{Sample $(\theta_e - \omega_R)$ profile for (a) large $\tau_1$, $\tau_3$  (b) small $\tau_1$, $\tau_3$}
    \label{fig:rate_profile}
\end{figure}

The the sample $(\theta_e - \omega_R)$ graph when $\alpha_R = 0.002$, $\tau_1 = 5$, $\tau_3 = 7$, and $\omega_{R_\text{max}} = 0.01745$ are demonstrated in Figure.~\ref{fig:rate_profile}(a). For comparison, the regulated rate profile (maximum profile) for a bang-bang structure without slew rate constraint is plotted as a blue solid line. If we reduce the time duration of linearly varying acceleration segments,e.g., $\tau_1 = 0.5$ and $\tau_3 = 0.5$, it converges the ideal bang-bang structure profile as shown in Figure.~\ref{fig:rate_profile}(b).

\textbf{Remark 5.} It is clearly shown that the trapezoidal profile is located inside the feasible region defined by a bang-bang profile with lower curvature. In addition, the duration for each segment is finite meaning that if the states follow this profile, $\theta_e$ will converge to 0 within finite time.

\textbf{Modified Trapezoidal Profile: } When we compute the $\dot{\boldsymbol{\omega}}_R$, some terms are function of $\dot{\hat{\textbf{\textit{e}}}}$, as will be shown in following subsection. If we use the profile converging faster than the linear $(\theta_e, \omega_R)$ profile, the $\dot{\boldsymbol{\omega}}_R$ term will chatter due to the larger value of $\dot{\hat{\textbf{\textit{e}}}}$, especially if the control frequency is low. Therefore, we replace the profile for the region $0 \leq \theta_e < \theta_1 $ as
\begin{equation} \label{eq:modified_profile}
    \omega_R = \begin{cases}
        \sqrt{\frac{\alpha_R}{\theta_1}}\theta_e, \quad \theta_{1} = \frac{1}{6}\alpha_R\tau_1^2 & \mbox{if } 0 \leq \theta_e < \theta_1 \\
        \text{Eq.~\eqref{eq:wr_profile} with } \omega_1 =  \sqrt{\frac{\alpha_R}{\theta_1}}\theta_1& \mbox{else}
    \end{cases}
\end{equation}
The the sample $(\theta_e - \omega_R)$ graph when $\alpha_R = 0.002$, $\tau_1 = 5$, $\tau_3 = 7$, and $\omega_{R_\text{max}} = 0.01745$ are demonstrated in Figure.~\ref{fig:rate_profile}(a). 

\textbf{Remark 6.} The profile is located inside both the bang-bang and trapezoidal profiles, with the first segment having a linear shape as expected. Note that this profile no longer satisfies finite-time convergence but only achieves exponential convergence.

%----------------------------------------------------------------------------------------
%	Rate Feedback Control
%----------------------------------------------------------------------------------------
\subsection{Rate Feedback Control and Implementation}
In order to compute the Eq.~\eqref{eq:control_law}, we now need a way to compute $\dot{\boldsymbol{\omega}}_R$. Based on the Eq.~\eqref{eq:wr_profile}, $\omega_R$ is a function of $\theta_e$, $\alpha_R$, and $\omega_{R_\text{max}}$, so the time derivative of it becomes:
\begin{equation} \label{eq:wr_dot}
    \dot{\boldsymbol{\omega}}_R = \left( \frac{\partial \omega_R}{\partial \theta_e} \dot \theta_e +  \frac{\partial \omega_R}{\partial \alpha_R} \dot \alpha_R +  \frac{\partial \omega_R}{\partial \omega_{R_\text{max}} } \dot{\omega}_{R_\text{max}}  \right) \hat{\textbf{\textit{e}}} + \omega_R \dot{\hat{\textbf{\textit{e}}}}
\end{equation}
Instead of computing the analytic expression, one can obtain the partial derivative using numerical differentiation as follows:
\begin{equation} \label{eq:wr_partial}
    \frac{\partial \omega_R}{\partial (\cdot)} \approx \frac{\omega_R \left( (\cdot) + \varepsilon \right) - \omega_R }{\varepsilon}
\end{equation}
This can be done efficiently since the original expression has a closed-form solution, and an accurate result can be obtained by choosing a small $\varepsilon \ll 1$, e.g., $\varepsilon \approx 1e^{-7}$. Such a small variation can be used without concern for noise amplification since they are partial derivatives at fixed instantaneous variables.

The only problem is the computation of $\dot{\alpha}_R$. We could obtain the analytic expression as:
\begin{equation} \label{eq:alpha_R_dot}
    \dot \alpha_R = \frac{\dot \theta_e }{\eta} (\alpha_\text{max} - \alpha_\text{min})  -\frac{\sigma \alpha_\text{max}}{\lVert J\hat{\textbf{\textit{e}}} \rVert^2} 
    \hat{\textbf{\textit{e}}}^\top J^\top J \dot{\hat{\textbf{\textit{e}}}} 
    - \frac{\sigma \gamma}{\lVert J\hat{\textbf{\textit{e}}} \rVert} \frac{d}{dt}\left(  
    \lVert J {}^\mathcal{B} \dot{\boldsymbol{\omega}}_\mathcal{D} \rVert + 
    \lVert \boldsymbol{\omega}_\mathcal{B} \times J \boldsymbol{\omega}_\mathcal{B} \rVert \right)
\end{equation}
but the last two terms require $\ddot {\boldsymbol{\omega}}_\mathcal{D}$ and $\dot{\boldsymbol{\omega}}_\mathcal{B}$ which we do not generally know. This paper numerically obtained the derivatives with the time step of the control frequency. Note that the numerical differentiation tends to be less sensitive due to the small $\sigma$ value, which helps mitigate noise amplification that could dominate the remaining terms.

With all the results the torque command satisfying the control constraint can be obtained as follows:
\begin{equation} \label{eq:control_trq}
    \textbf{\textit{u}}_\text{cmd} = \begin{cases}
        \textbf{\textit{u}} & \mbox{if } u \leq u_\text{max} \\
        u_\text{max} \hat{\textbf{\textit{u}}} & \mbox{if } u > u_\text{max} \\
    \end{cases}
\end{equation}
In addition, following table explains the overall torque command generation sequence

\begin{algorithm} 
    \caption{Torque Command Generation}
  \begin{algorithmic}[1]
    \REQUIRE Satellite Parameters: $(J, \omega_\text{max}, u_\text{max})$, 
    Control Parameters: $(d_\text{max}, \gamma, \eta, \beta_1, \beta_2, \tau_1, \tau_3)$
    \INPUT Satellite States: $(\bar{\textbf{\textit{q}}}_\mathcal{B}, \boldsymbol{\omega}_\mathcal{B})$, 
    Desired Frame Profiles: $(\bar{\textbf{\textit{q}}}_\mathcal{D}, \boldsymbol{\omega}_\mathcal{D}, \dot{\boldsymbol{\omega}}_\mathcal{D})$
    \OUTPUT Command Torque: $(\textbf{\textit{u}}_\text{cmd})$
    \STATE Compute $\bar{\textbf{\textit{q}}}_e$, $\boldsymbol{\omega}_e$ using Eq.~\eqref{eq:error_state}.
    \STATE Compute $\hat{\textbf{\textit{e}}}$, $\theta_e$ using Eq.~\eqref{eq:eigenaxis} and $\dot{\hat{\textbf{\textit{e}}}}$, $\dot{\theta}_e$ using Eqs.~\eqref{eq:errorangle_dot}, ~\eqref{eq:eigenaxis_dot}
    \STATE Compute $\alpha_\text{max}$, $\alpha_R$, and $\dot{\alpha}_R$ using Eqs.~\eqref{eq:alpha_max},~\eqref{eq:alpha_R}, and ~\eqref{eq:alpha_R_dot}.
    \STATE Compute $\omega_{R_\text{max}}$ and $\dot{\omega}_{R_\text{max}}$ using Eqs.~\eqref{eq:omega_max},~\eqref{eq:omega_max_dot}
    \STATE Compute $\omega_R$ either using the trapezoidal profile Eq.~\eqref{eq:wr_profile} or modified trapezoidal profile Eq.~\eqref{eq:modified_profile}
    \STATE Compute $\partial_{(\cdot)}\omega_R$ numerically as Eq.~\eqref{eq:wr_partial} and compute $\dot{\boldsymbol{\omega}}_R$ as Eq.~\eqref{eq:wr_dot}
    \STATE Compute control law $\textbf{\textit{u}}$ as Eq.\eqref{eq:control_law} and apply control torque $\textbf{\textit{u}}_\text{cmd}$ as Eq.\eqref{eq:control_trq}.
  \end{algorithmic}
\end{algorithm}

%----------------------------------------------------------------------------------------
%	Simulation Result
%----------------------------------------------------------------------------------------
\section{Simulation Result}
\subsection{Simulation Environments}
The following parameters are used for the simulation of the entire scenario. Note that the satellite parameters are adopted from the paper[], and reduced values are used for both the maximum slew rate and torque level to better demonstrate the effectiveness of algorithm.
\begin{table}[h]
    \centering
    \begin{tabular}{c|c|c}
    Parameters & Value  & Description \\ \hline
    $J$ & $\begin{bmatrix}
        21400 & 2100 & 1800 \\ 2100 & 20100 & 500 \\ 1800 & 500 & 5000 
    \end{bmatrix}$  & Moment of inertia of satellite $[\text{kg/m}^2]$ \\
    $\omega_\text{max}$ & 3 & Maximum Slew Rate $[\text{deg/sec}]$ \\
    $u_\text{max}$ & 150 & Maximum control torque level  $[\text{N}\cdot\text{m}]$ \\
    $d_\text{max}$ & 2 & Maximum disturbance torque level $[\text{N}\cdot\text{m}]$ \\
    $\gamma$ & 0.99 & Reduced torque ratio $[-]$ \\
    $\eta$ & 0.05 & Acceleration transition threshold $[\text{deg}]$ \\
    $\beta_1$ & 2 & Sliding control gain $[-]$ \\
    $\beta_2$ & 0.5 & Sliding control exponent $[-]$ \\
    $\tau_1$ & 1 & Duration of the first linear segment $[\text{sec}]$ \\
    $\tau_3$ & 1 & Duration of the last linear segment $[\text{sec}]$ \\
    \hline
    \end{tabular}
    \caption{Parameters for the direct resupply method analysis}
\end{table}
The control is updated at a 10 Hz frequency, and the dynamics are integrated using the RK4 method with a time step of 0.01 seconds. Lastly, the following disturbance, having a near-orbital frequency, is considered for the entire simulation.
\begin{equation}
    \textbf{\textit{d}}^\mathcal{B}(t) = \begin{bmatrix}
        1.1 \sin\left(0.0012t + \frac{\pi}{6}\right) \\
        0.9 \sin\left(0.0010t\right) \\
        1.0 \sin\left(0.0013t + \frac{\pi}{2}\right)
    \end{bmatrix}
\end{equation}

\subsection{Three-axis Rest-to-Rest Maneuver}
This simulation is designed to validate the performance of the control law for a rest-to-rest maneuver, i.e., initial and final angular velocities are zero. The first set of figures presents the results when the trapezoidal profile is used for a 90-degree roll-axis maneuver. Each figure shows the error angle ($\theta_e$) and the norm of the error rate ($\lVert \boldsymbol{\omega}_e\rVert$), the norm of the body angular velocity ($\lVert \boldsymbol{\omega}_\mathcal{B} \rVert$), the torque command for each axis ($u_i,\ i=x,y,z$), and the norm of the torque command ($\lVert \textbf{\textit{u}} \rVert$). The next set of figures illustrates the results for the modified trapezoidal profile under the same scenario. The results indicate that both profiles satisfy the slew-rate and control constraints, but the original trapezoidal profile exhibits larger chattering in the torque command due to previously explained reasons.
\begin{figure}[H]
    \centering
    \includegraphics[width=0.42\textwidth]{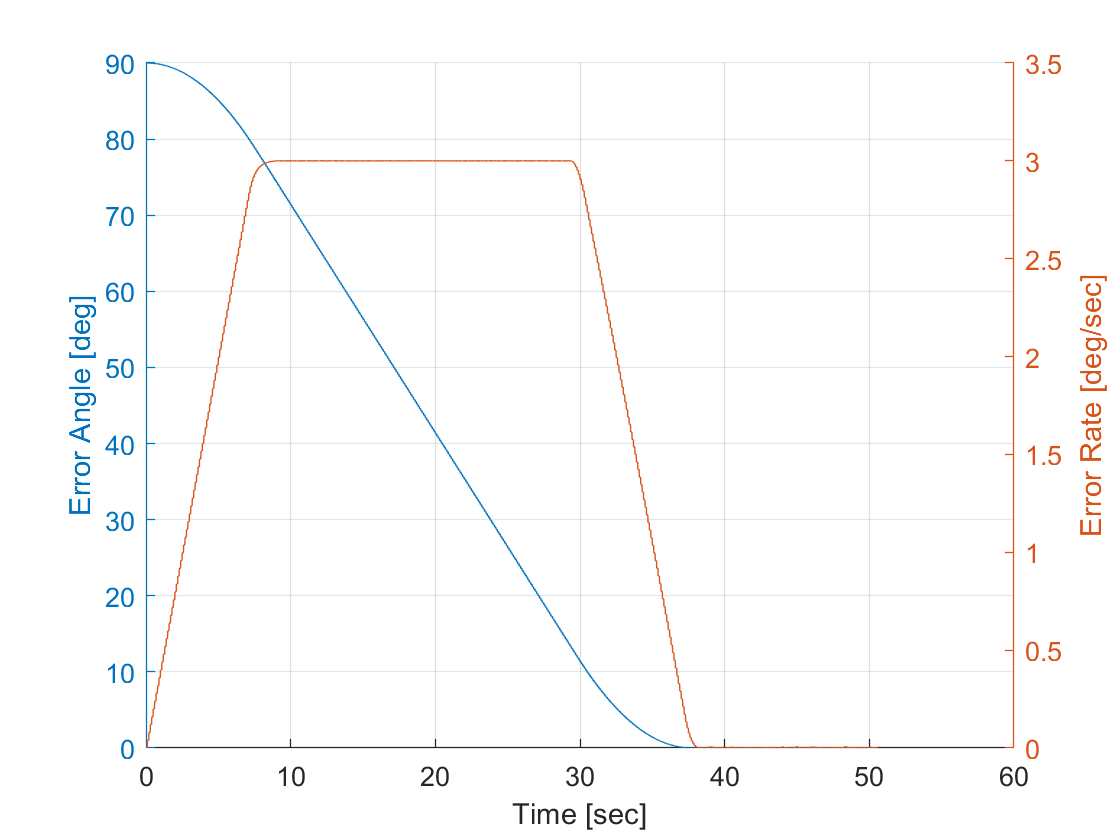}
    \includegraphics[width=0.42\textwidth]{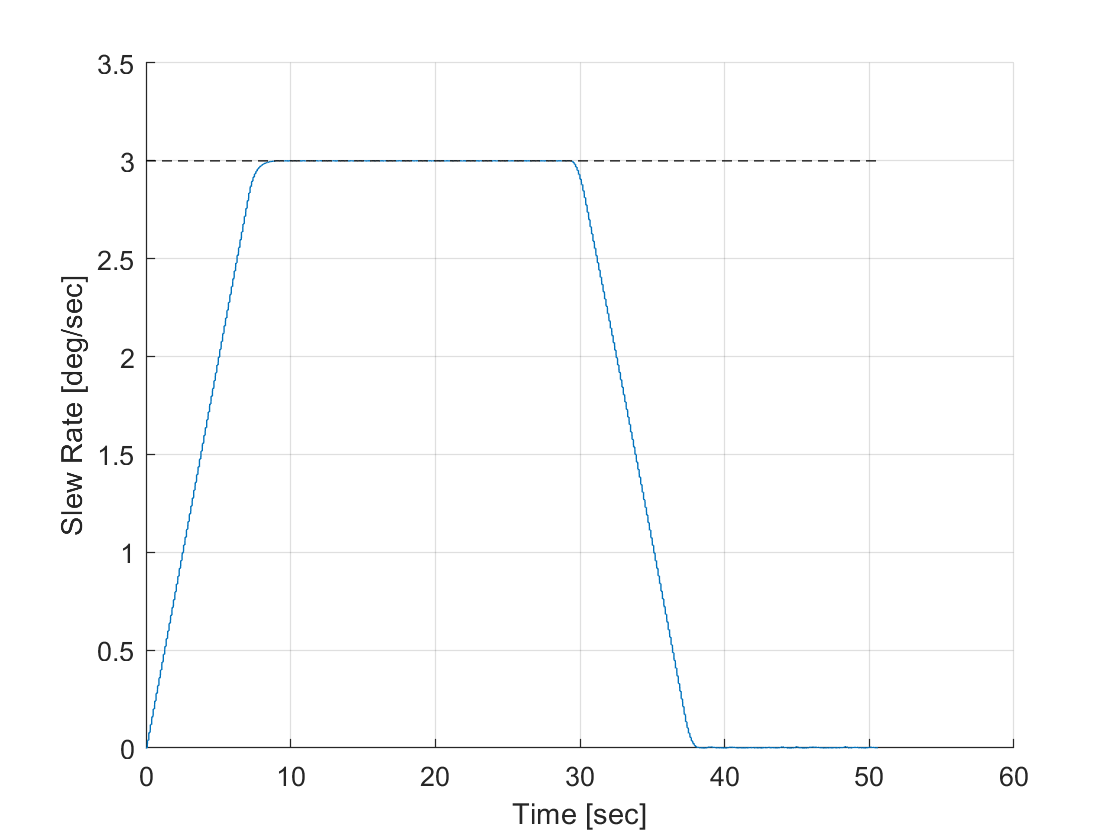}
    \includegraphics[width=0.42\textwidth]{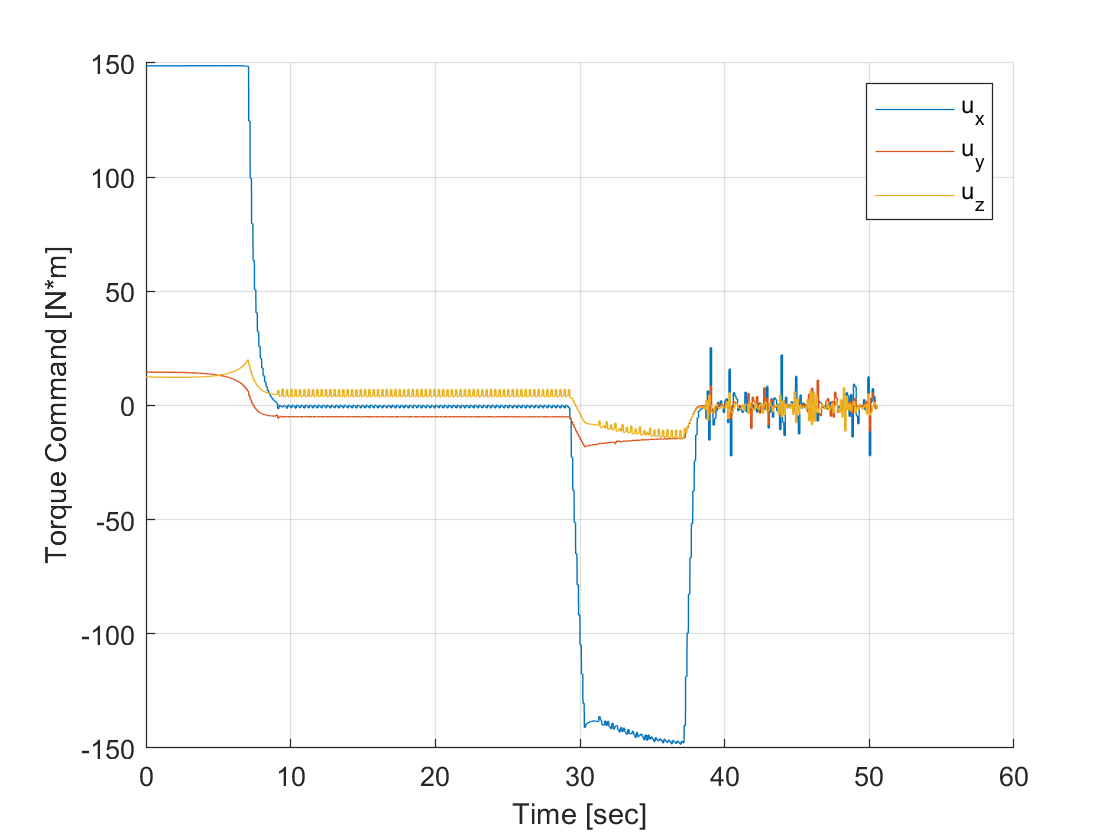}
    \includegraphics[width=0.42\textwidth]{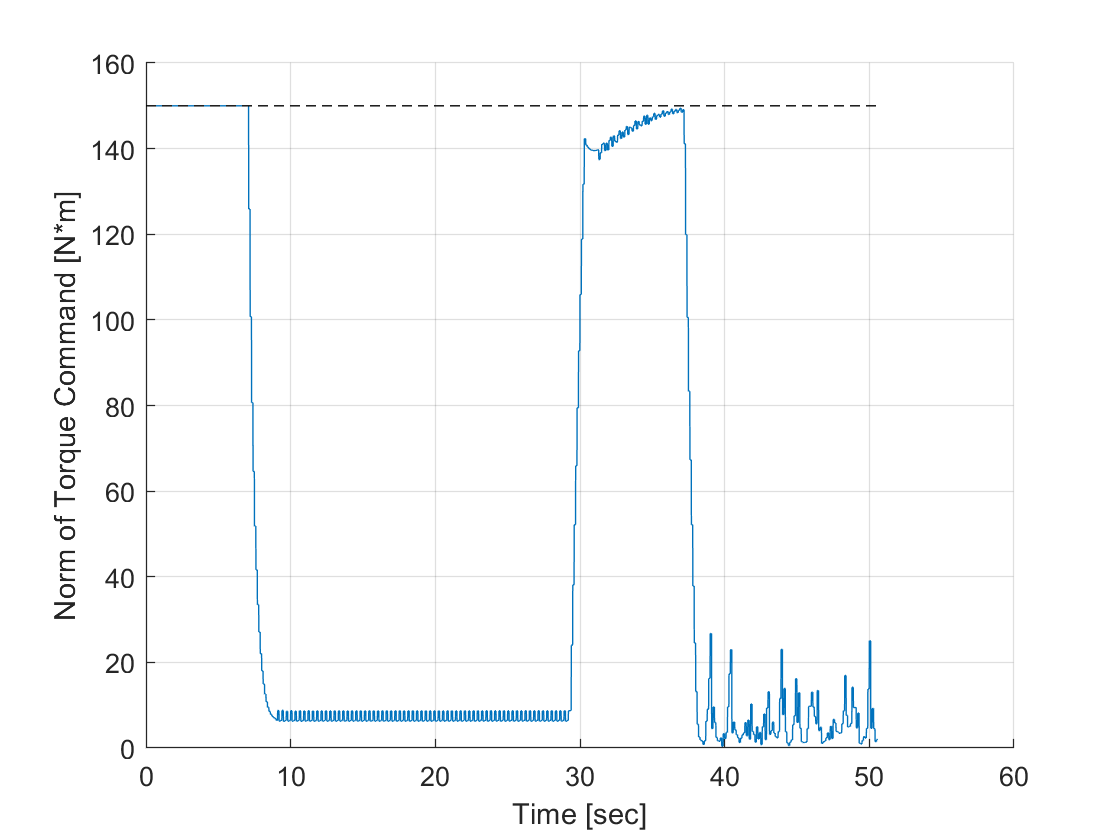}
    \caption{Result of Trapezoidal Profile (a) $\theta_e$, $\lVert \boldsymbol{\omega}_e\rVert$  (b) $\lVert \boldsymbol{\omega}_\mathcal{B} \rVert$ (c) $u_i,\ i=x,y,z$ (d) $\lVert \textbf{\textit{u}} \rVert$}
    \label{fig:r2r_trape}
\end{figure}
\begin{figure}[h]
    \centering
    \includegraphics[width=0.42\textwidth]{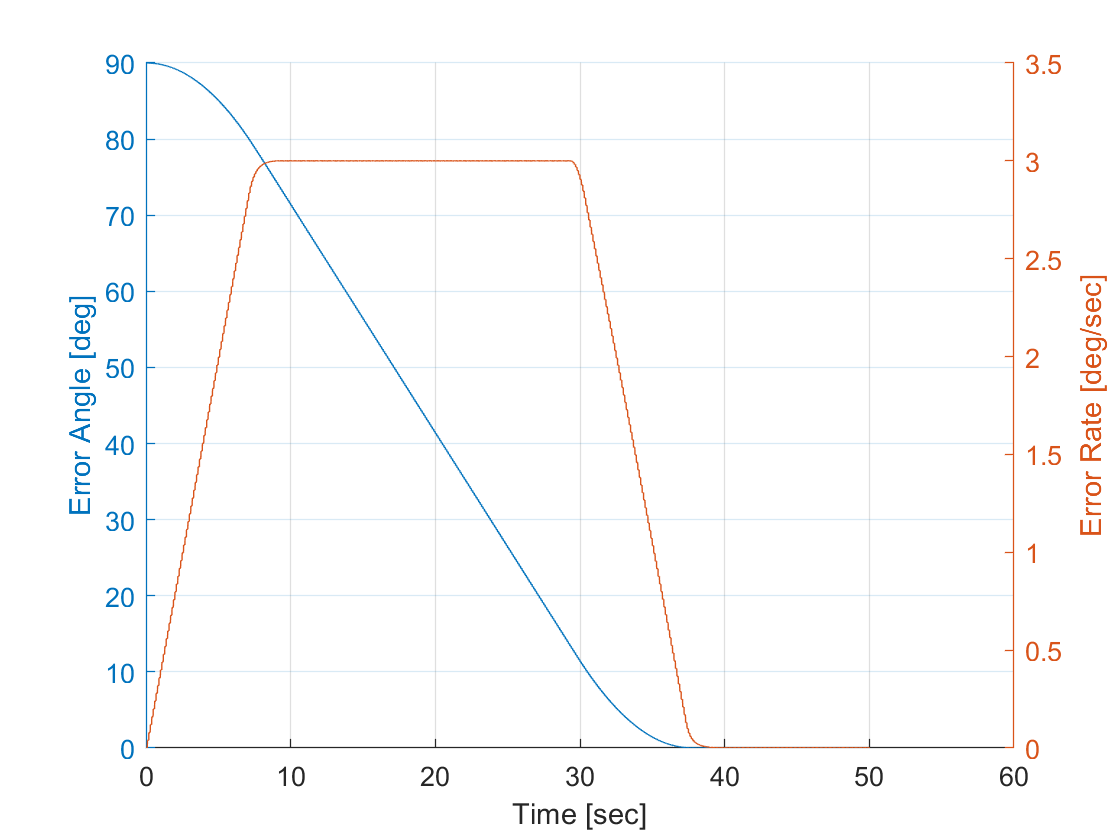}
    \includegraphics[width=0.42\textwidth]{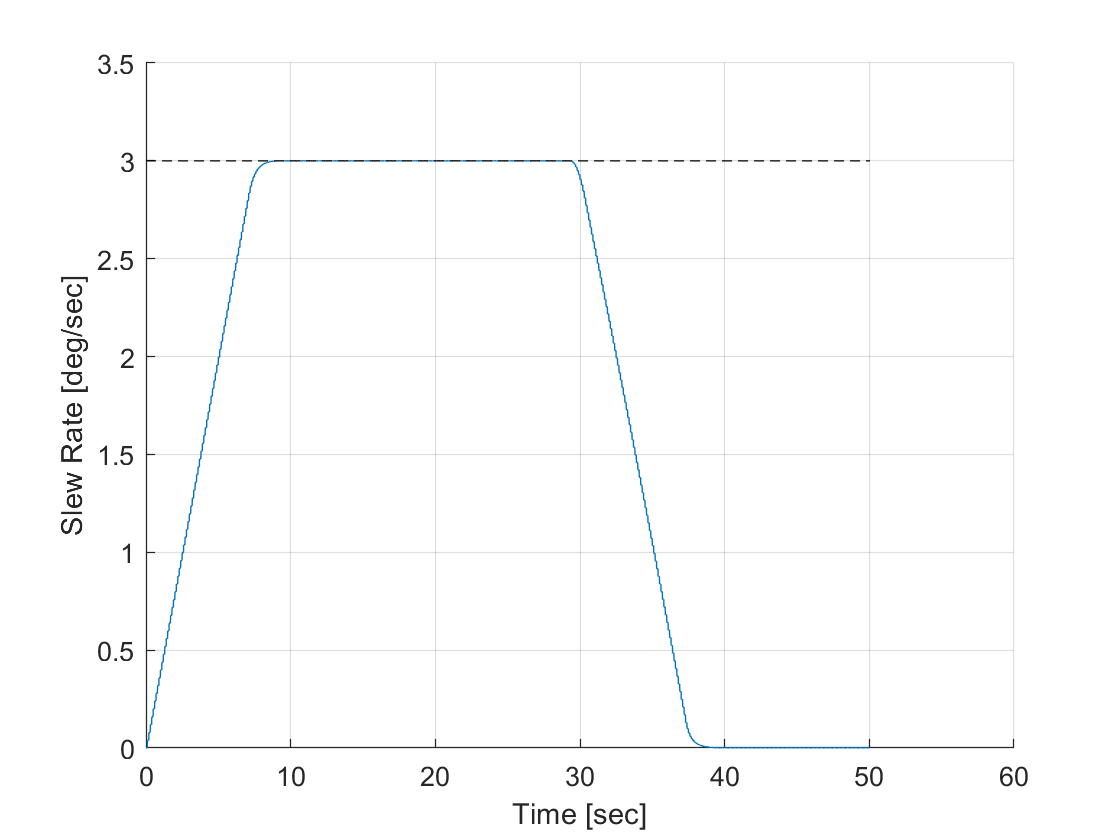}
    \includegraphics[width=0.42\textwidth]{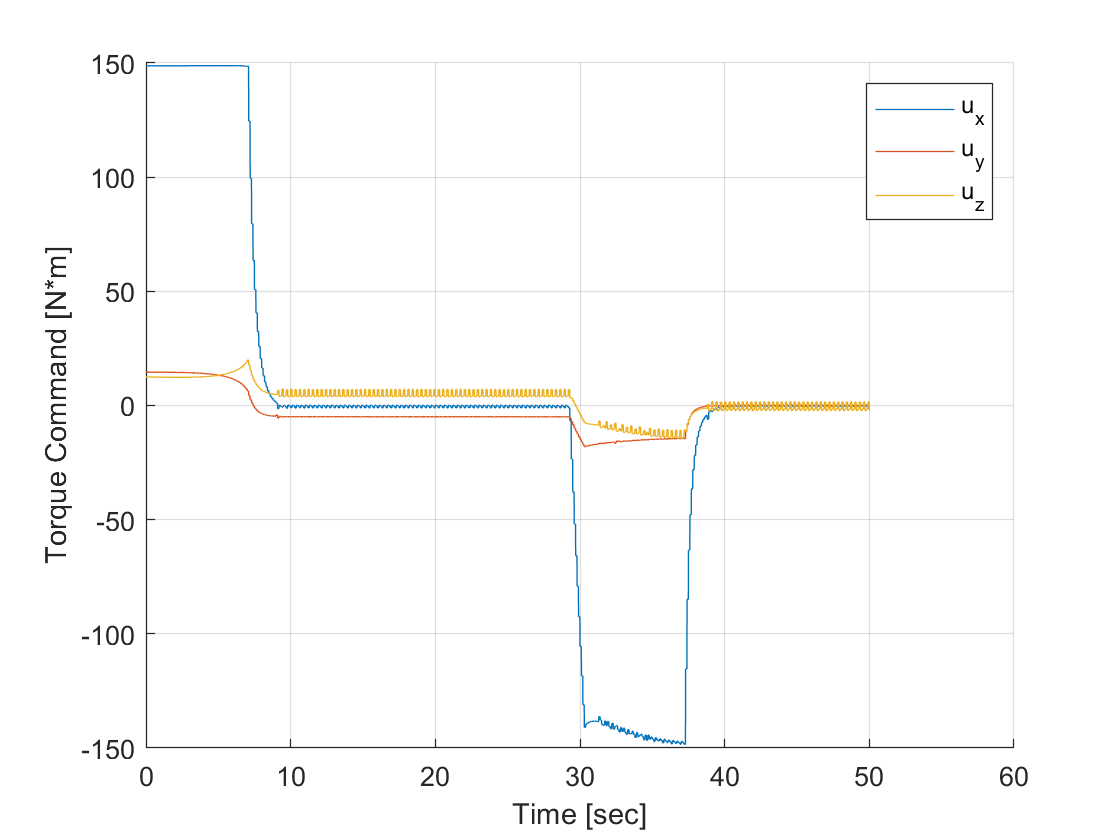}
    \includegraphics[width=0.42\textwidth]{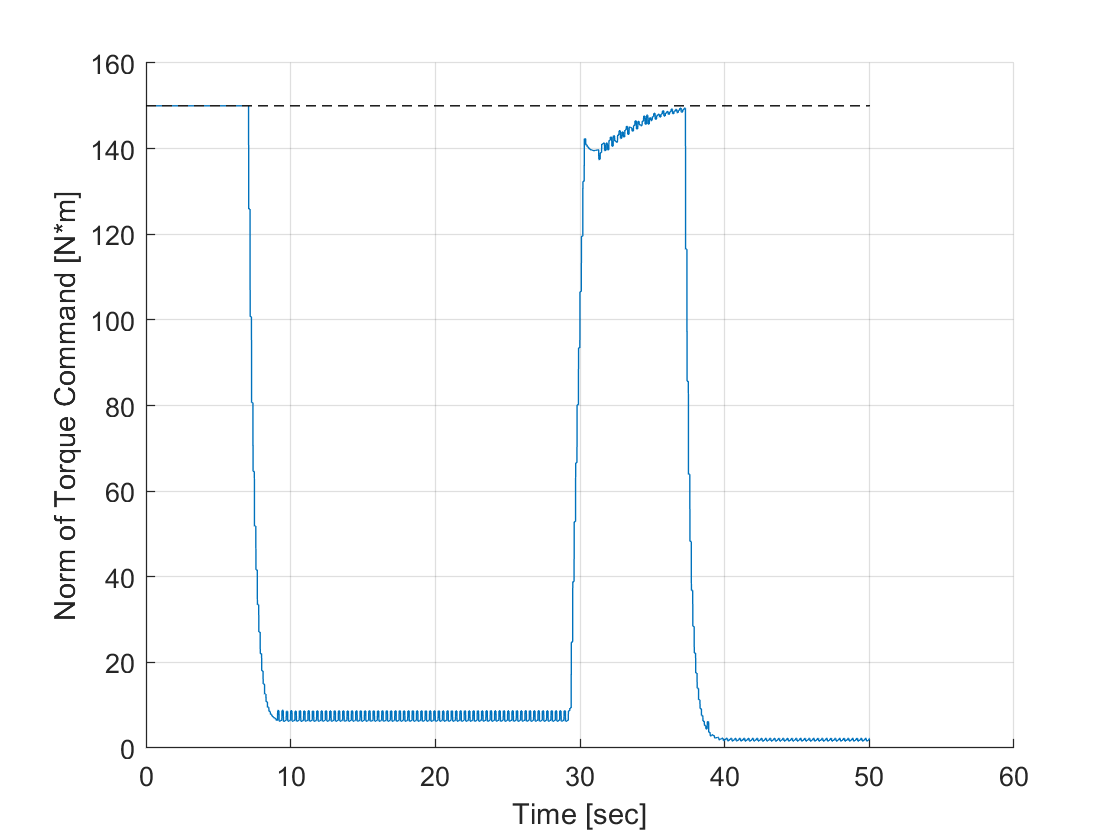}
    \caption{Result of Modified Trapezoidal Profile (a) $\theta_e$, $\lVert \boldsymbol{\omega}_e\rVert$  (b) $\lVert \boldsymbol{\omega}_\mathcal{B} \rVert$ (c) $u_i,\ i=x,y,z$ (d) $\lVert \textbf{\textit{u}} \rVert$}
    \label{fig:r2r_mod}
\end{figure}

Note that if we increase the control frequency, then the chattering issue of the trapezoidal profile will be diminished. That is the profile works well in a continuous manner, it is not feasible in practice.
\begin{figure}[H]
    \centering
    \includegraphics[width=0.42\textwidth]{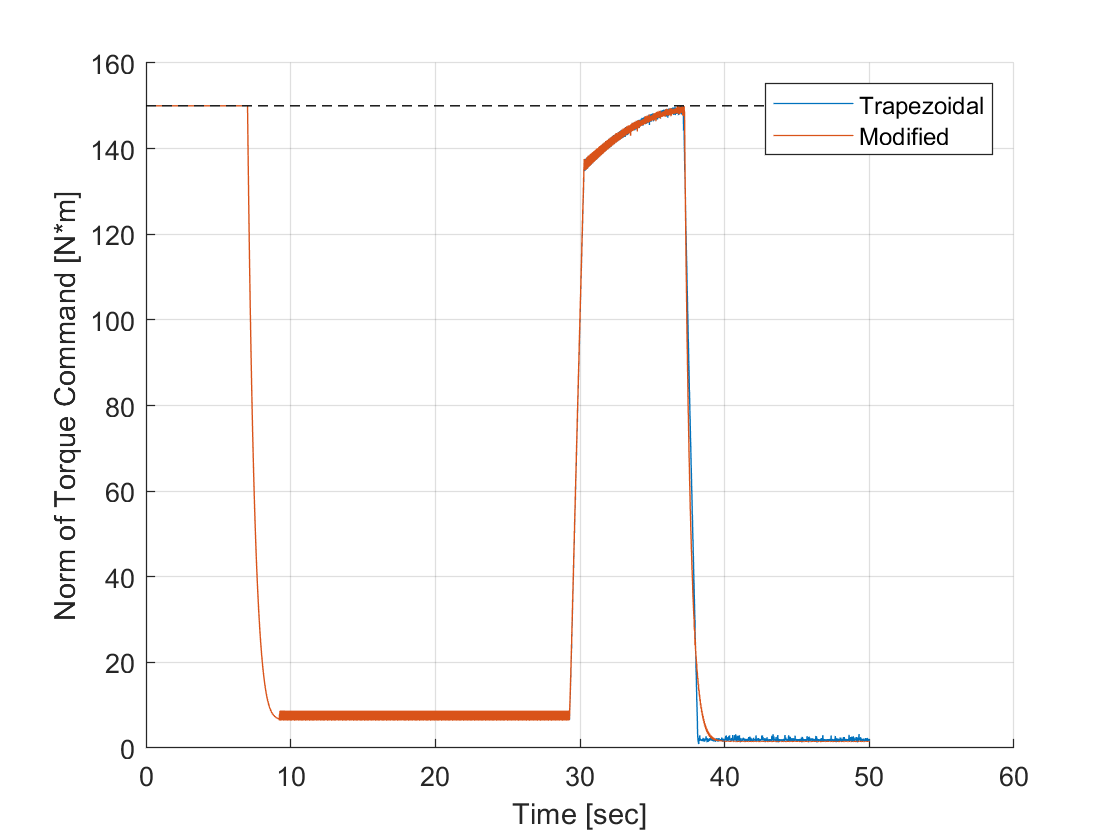}
    \caption{Norm of torque command comparison when control frequency is 100 Hz}
\end{figure}

The next figures show the time required to reorient the satellite for different maneuver angles and rotational axes. The reorientation time is measured until the error angle becomes less than $0.01\ \text{deg}$ and the error rate becomes less than $0.01\ \text{deg/sec}$, which are practical values. As can be seen, the reorientation time is almost identical for both profiles, with the modified trapezoidal profile taking at most 0.6 seconds longer. However, the modified trapezoidal profile exhibits a preferable torque profile after convergence, so will be used in the subsequent scenario.
\begin{figure}[H]
    \centering
    \includegraphics[width=0.42\textwidth]{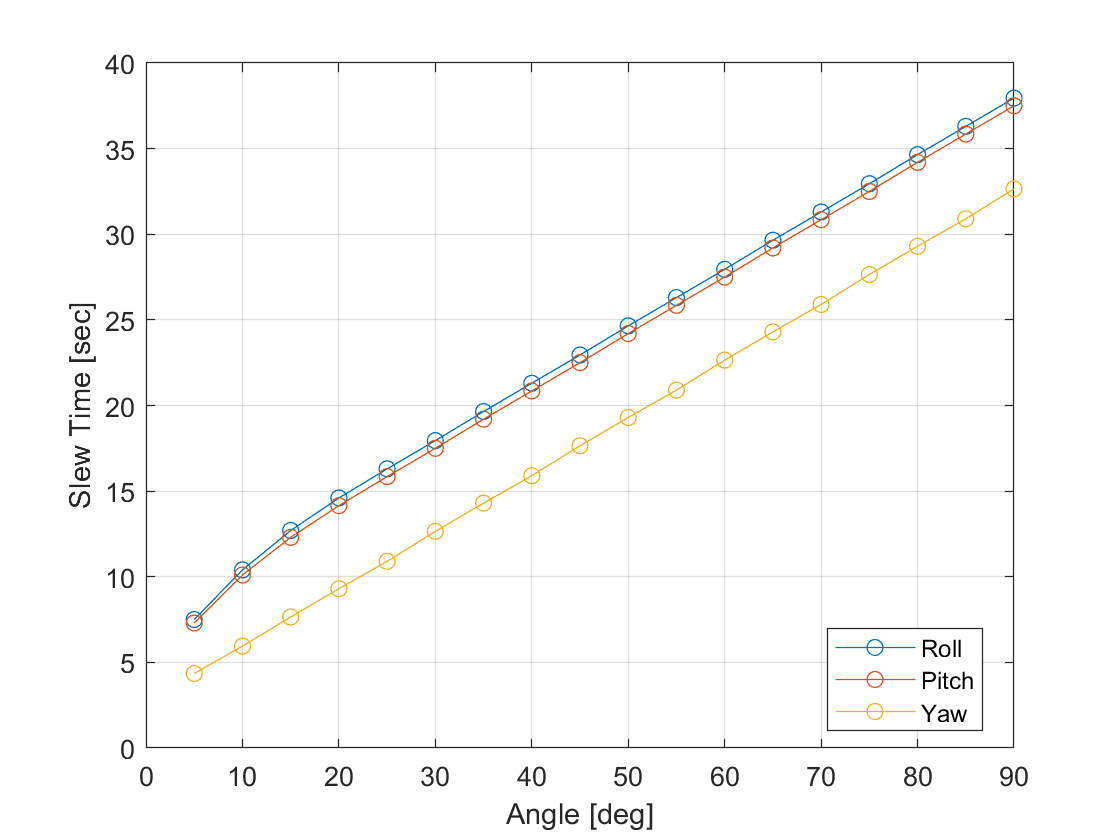}
    \includegraphics[width=0.42\textwidth]{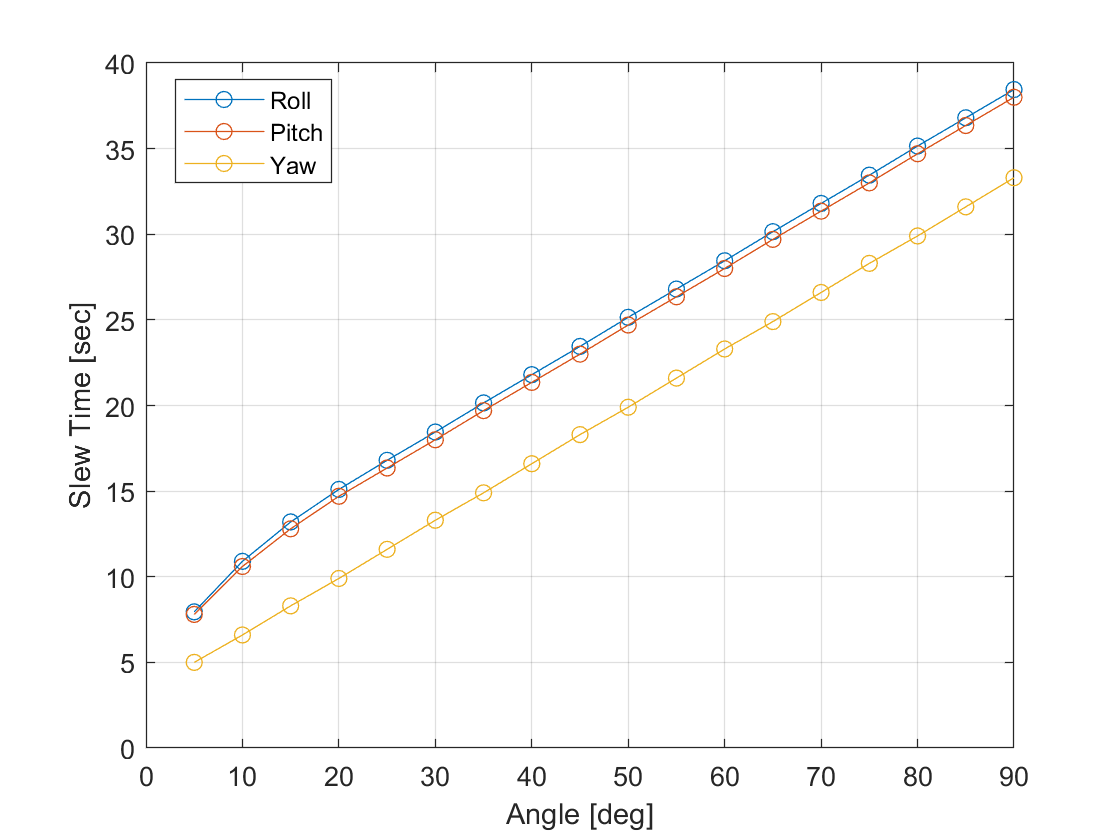}
    \caption{Slew Time for Each axis Maneuver (a) Trapezoidal Profile (b) Modified Trapezoidal Profile}
\end{figure}

\subsection{Two Successive SPOT Imaging}
This simulation is designed to verify the tracking performance of the control law for practical imaging operations. The scenario consists of two spotlight (SPOT) imaging operations (staring at a ground-fixed target), and the analytical expressions for $\bar{\textbf{\textit{q}}}_\mathcal{D}$, $\boldsymbol{\omega}_\mathcal{D}$, and $\dot{\boldsymbol{\omega}}_\mathcal{D}$ can be found in the referenced paper\cite{han2022analytical}. The altitude of the satellite is set to $500\ \text{km}$, and the look angles of the targets are set to $5$ and $25\ \text{deg}$ degrees, respectively, when the squint angle is $90\ \text{deg}$ degrees, as illustrated in the following figure.
\begin{figure}[H]
    \centering
    \includegraphics[width=0.42\textwidth]{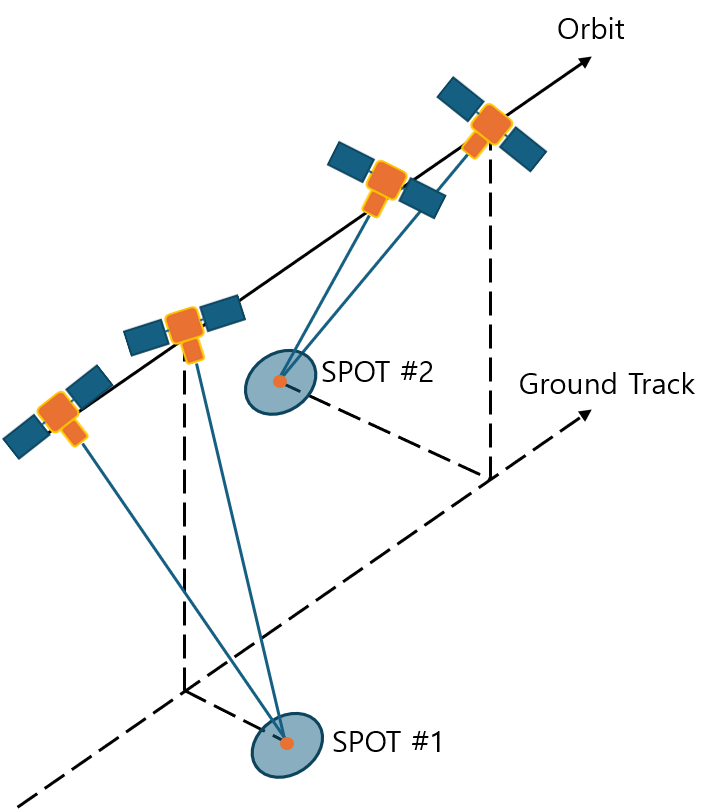}
    \caption{Two successive SPOT imaging operation}
\end{figure}

Based on the result, it is clear that the control law time-efficiently reorients the satellite while satisfying all the constraint requirements, even when the desired frame is time-varying. As shown in the rate profile graph, the angular velocity of $\mathcal{D}$ is around $0.8\ \text{deg/sec}$ which is significant, still the proposed method can accurately track this practical imaging operation. Note that the first figure shows the 3-2-1 Euler angles between $\mathcal{D}$ and $\mathcal{O}$, and $\mathcal{B}$ and $\mathcal{O}$ where $\mathcal{O}$ is the orbital frame.
\begin{figure}[H]
    \centering
    \includegraphics[width=0.42\textwidth]{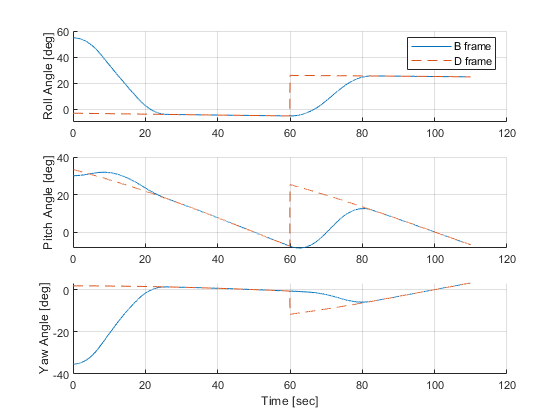}
    \includegraphics[width=0.42\textwidth]{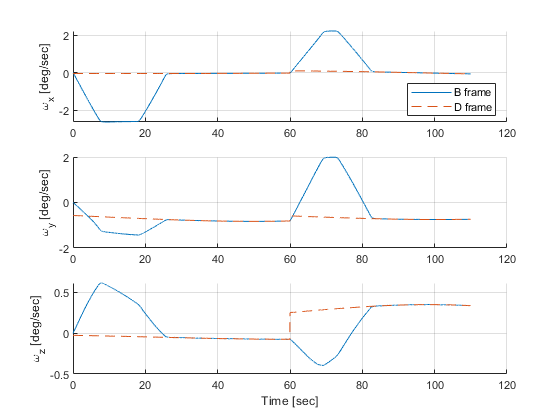}
    \includegraphics[width=0.42\textwidth]{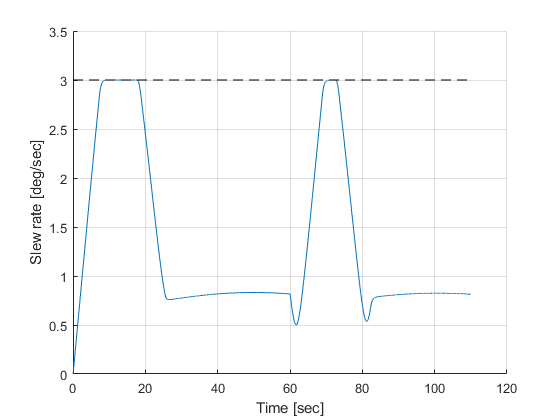}
    \includegraphics[width=0.42\textwidth]{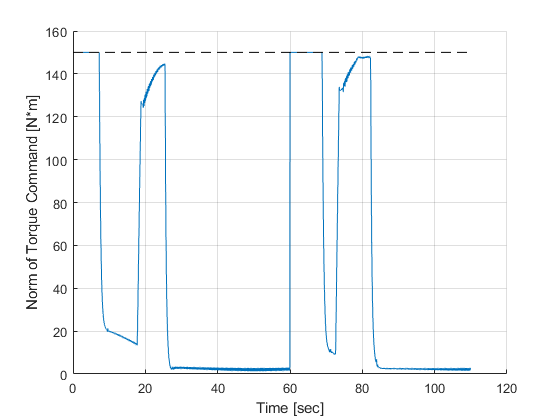}
    \caption{Successive SPOT Results (a) Euler 3-2-1 Angles, (b) $\boldsymbol{\omega}_\mathcal{B}$, $\boldsymbol{\omega}_\mathcal{D}$ (c) $\lVert \boldsymbol{\omega}_\mathcal{B} \rVert$ (d) $\lVert \textbf{\textit{u}} \rVert$}
\end{figure}

\subsection{Two Successive STRIP Imaging}
This simulation is designed to verify the tracking performance of the control law for special imaging operations. The scenario consists of two strip imaging operations (scanning the ground target profile), and the analytical expressions for $\bar{\textbf{\textit{q}}}_\mathcal{D}$, $\boldsymbol{\omega}_\mathcal{D}$, and $\dot{\boldsymbol{\omega}}_\mathcal{D}$ can be found in the referenced paper\cite{han2024strip}. The scenario is illustrated in the following figure, and it requires high tracking velocity when the ground scanning direction is opposite to the orbit direction. Note that a spherical Earth model is used.
\begin{figure}[H]
    \centering
    \includegraphics[width=0.35\textwidth]{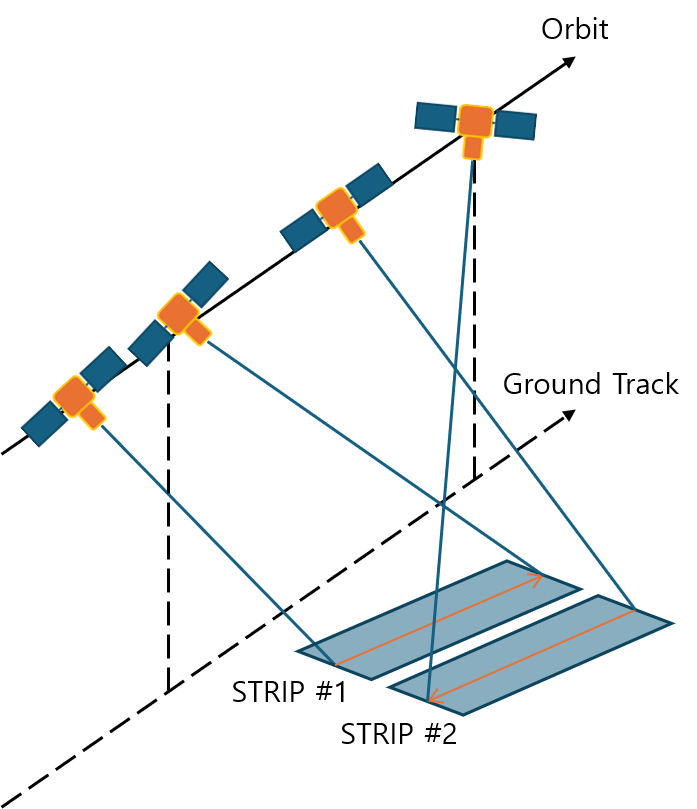}
    \caption{Two successive STRIP imaging operation}
\end{figure}

Based on the results, it is clear that the control law time-efficiently reorients the satellite while satisfying all constraint requirements, even when the desired frame is rapidly time-varying. Although this kind of operation is not common for many EOS, some agile EOS require this capability, and this control law is capable of supporting such imaging operations.
\begin{figure}[h]
    \centering
    \includegraphics[width=0.42\textwidth]{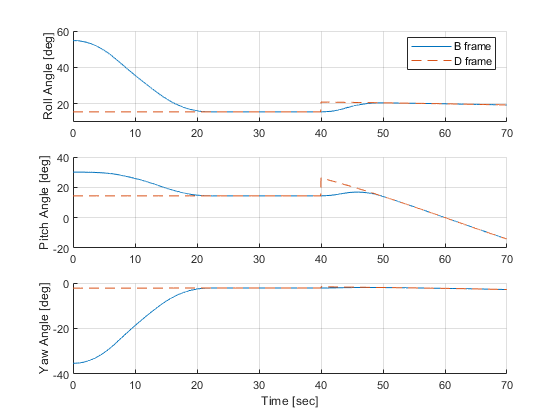}
    \includegraphics[width=0.42\textwidth]{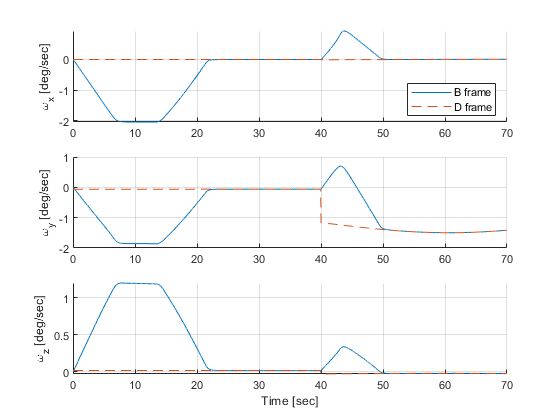}
    \includegraphics[width=0.42\textwidth]{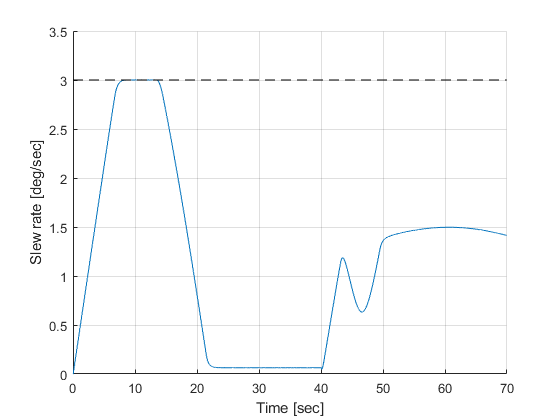}
    \includegraphics[width=0.42\textwidth]{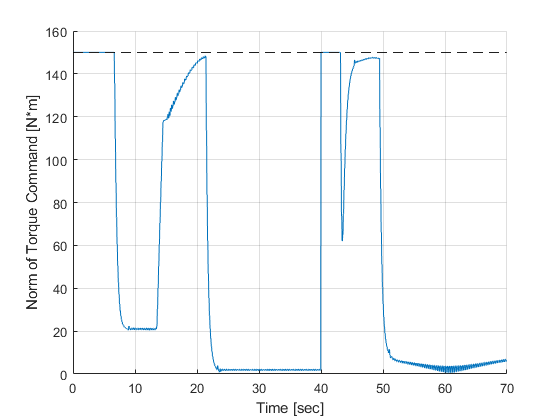}
    \caption{Successive STRIP Results (a) Euler 3-2-1 Angles, (b) $\boldsymbol{\omega}_\mathcal{B}$, $\boldsymbol{\omega}_\mathcal{D}$ (c) $\lVert \boldsymbol{\omega}_\mathcal{B} \rVert$ (d) $\lVert \textbf{\textit{u}} \rVert$}
\end{figure}
%----------------------------------------------------------------------------------------
%	Conclusion
%----------------------------------------------------------------------------------------
\section{Conclusion}
In this paper, a practical time-efficient attitude tracking control algorithm is proposed, capable of handling both slew-rate and control constraints. The regulating rate is designed to achieve time-efficient maneuvers while satisfying the slew-rate and control torque constraints. The performance and stability of the proposed control law are demonstrated through a practical imaging operation scenario. Based on the simulation results, the proposed law is capable of effectively supporting various imaging operations.

Future work will involve adopting an improved sliding mode control structure to reduce the chattering problem, and testing the performance of the control law with practical sensor noise.

\section{Appendix} 
\subsection{Derivation of Eq.(28)}
\subsubsection{First segment}
The duration of the first segment ($\tau_1$) and the maximum acceleration ($\alpha_{R}$) are given. Let $\tau = t - t_0 = t$ be the elapsed time since the beginning of the first segment.  Then, the acceleration, velocity, and position profile of the first segment becomes:
\begin{equation}
\begin{aligned}
    \alpha(\tau) &= \alpha_{R}\frac{\tau}{\tau_1} \\
    \omega(\tau) &= \omega(0) + \alpha_{R}\frac{\tau^2}{2\tau_1} = \alpha_{R}\frac{\tau^2}{2\tau_1}\\
    \theta(\tau) &= \theta(0) + \alpha_{R}\frac{\tau^3}{6\tau_1} = \alpha_{R}\frac{\tau^3}{6\tau_1}
\end{aligned}
\end{equation}
since $\omega(0) = \theta(0) = 0$ by the definition. Substituting $\tau =\tau_1$ gives the acceleration, velocity, and position at the end of the first segment as:
\begin{equation}
\begin{aligned}
    \alpha_1 = \alpha_{R},\quad
    \omega_1 = \alpha_{R}\frac{\tau_1}{2},\quad
    \theta_1 = \alpha_{R}\frac{\tau_1^2}{6}
\end{aligned}
\end{equation}
Therefore, the range of $\theta(t)$ for the first segment becomes $0\leq \theta \leq \theta_1$, and for given $\theta(\tau)$:
\begin{equation}
    \theta(\tau) \implies \tau = \left( \frac{6\theta(t)\tau_1}{\alpha_{R_\text{max}}}\right)^{\frac{1}{3}} \implies \omega(\tau)
\end{equation}
$\omega(\tau)$ can be obtained.

\subsubsection{Second segment}
Since we know the duration of the third segment ($\tau_3$), and the rate at the end of the third segment ($\omega_3 = \omega_{R_\text{max}}$), the rate at the end of the second segment (or the beginning of the third segment) is
\begin{equation}
    \omega_2 = \omega_{R_\text{max}} - \frac{1}{2}\alpha_{R}\tau_3
\end{equation}
Then, the duration of the second segment ($\tau_2$) becomes:
\begin{equation}
    \tau_2 = \frac{\omega_2 - \omega_1}{\alpha_R}
\end{equation}
since the acceleration is constant during the second segment. Let $\tau = t - t_1 = t - \tau_1$, i.e., the elapsed time from the beginning of the second segment. The acceleration, velocity, and position profile of the second segment become:
\begin{equation}
\begin{aligned}
    \alpha(\tau) &= \alpha_{R} \\
    \omega(\tau) &= \omega_1 + \alpha_{R}\tau \\
    \theta(\tau) &= \theta_1 + \omega_1 \tau + \frac{1}{2}\alpha_{R}\tau^2
\end{aligned}
\end{equation}
Therefore, the range of $\theta(\tau)$ for the second segment becomes $\theta_1 \leq \theta \leq \theta_2$ where 
\begin{equation}
    \theta_2 = \theta_1 + \omega_1 \tau_2 + \frac{1}{2}\alpha_R \tau_2^2
\end{equation}
and for given $\theta(\tau)$:
\begin{equation}
    \theta(\tau) \implies \tau = \frac{-\omega_1 + \sqrt{\omega_1^2 + 2\alpha_R(\theta(t) - \theta_1)}}{\alpha_R}
    \implies \omega(\tau)
\end{equation}
$\omega(\tau)$ can be obtained.

\subsubsection{Third segment} Let $\tau = t_3 - t = \tau_1 +  \tau_2 + \tau_3 - t$,  i.e., the remaining time to the end of the third segment. The reason for this backward time integration is to directly make $\theta(t)$ have a cubic polynomial expression in the form $x^3 + px + q = 0$. The acceleration, velocity, and position profile of the third segment become:
\begin{equation}
\begin{aligned}
    \alpha(\tau) &= \alpha_{R}\frac{\tau}{\tau_3} =  \alpha_{R}\frac{t_3 - t}{\tau_3}\\
    \omega(\tau) &= \omega_3 - \frac{\alpha_{R}}{2\tau_3}\tau^2 \\
    \theta(\tau) &= \theta_3 - \omega_3 \tau + \frac{\alpha_{R}}{6\tau_3}\tau^3
\end{aligned}
\end{equation}
Note that the following relation is used for the integration of velocity and position:
\begin{equation}
    \int_{t}^{t_3} \dot \omega\ dt = \omega_3 - \omega(t) = \int_{t}^{t_3} \alpha_R(t_3 - s)\ ds = - \int_{\tau}^0 \alpha_R \sigma\ d\sigma = \frac{1}{2}\alpha_R \tau^2
\end{equation}
\begin{equation}
    \int_{t}^{t_3} \dot \theta\ dt = \theta_3 - \theta(t) = \int_{t}^{t_3} \omega_3 - \frac{1}{2}\alpha_R(t_3-s)^2\ ds = - \int_{\tau}^0 \omega_3 - \frac{1}{2}\alpha_R\sigma^2\ d\sigma 
    = \omega_3\tau - \frac{1}{6}\alpha_R \tau^3
\end{equation}
Applying $\tau = \tau_3$ gives the $\theta_3$ since $\tau$ is the remaining time to $t_3$.
\begin{equation}
\begin{aligned}
    \theta_2 = \theta(\tau_3) = \theta_3 - \omega_3 \tau_3 + \frac{1}{6}\alpha_R \tau_3^2 \implies \theta_3 &= \theta_2 + \omega_3 \tau_3 - \frac{1}{6}\alpha_R \tau_3^2 \\
    &= \theta_2 + \left(\omega_2 + \frac{\alpha_R}{2\tau_3}\tau_3^2\right) \tau_3 - \frac{1}{6}\alpha_R \tau_3^2 \\
    &=  \theta_2 + \omega_2 \tau_3 + \frac{1}{3}\alpha_R \tau_3^2
\end{aligned}
\end{equation}
For given $\theta(t) \in [\theta_2,\ \theta_3]$:
\begin{equation}
    \tau^3 -\frac{6\omega_3 \tau_3}{\alpha_R}\tau + \frac{6\tau_3(\theta_3 - \theta)}{\alpha_R} = 0 \implies
    \tau^3 + p\tau + q =0,\ p = -\frac{6\omega_3 \tau_3}{\alpha_R},\ q = \frac{6\tau_3(\theta_3 - \theta)}{\alpha_R}
\end{equation}
Without proof, we assume $4p^3 + 27q^2 < 0$ and $\tau$ has one negative and two positive solutions. Then, the proper solution $\tau$ for a given $\theta(\tau)$ can be obtained using the trigonometric function\cite{korn2000mathematical}:
\begin{equation}
\theta(\tau) \implies
    \tau = 2\sqrt{-\frac{p}{3}}\cos\left(\frac{1}{3}\cos^{-1}\left(\frac{3q}{2p} \sqrt{\frac{-3}{p}}\right) -\frac{2\pi}{3}\right) \implies \omega(\tau)
\end{equation}
and $\omega(\tau)$ can be computed.

\subsubsection{Fourth segment} There is no more acceleration and the rate reaches its maximum value during the fourth segment, so for $\theta(t) > \theta_3$
\begin{equation}
\begin{aligned}
    \alpha(\tau) &= 0\\
    \omega(\tau) &= \omega_3 = \omega_{R_\text{max}}\\
\end{aligned}
\end{equation}

\bibliographystyle{AAS_publication}   % Number the references.
\bibliography{references}   % Use references.bib to resolve the labels.
\end{document}